%% file: Main.tex
\documentclass[12pt]{amsart}
\pdfoutput=1 
\usepackage[margin=1.25in]{geometry}
\usepackage{tangle}
\usepackage{amsmath,amssymb,amsthm,mathabx,mathtools}
\usepackage{mathrsfs}
\usepackage{tikz-cd}
\usepackage{adjustbox}
\usepackage{graphicx}
\usepackage[all,cmtip]{xy}
\usepackage[backref=page, bookmarks=false, hidelinks]{hyperref}
\hypersetup{
    pdfauthor={Angel Toledo},
    pdftitle={Davydov-Yetter cohomology for tensor triangulated categories},
    colorlinks=true,
    citecolor=blue,
    linkcolor=black,
}
\newcommand{\hzero}{H^{0}}
\newcommand{\cc}{\mathscr{C}(k)}
\newcommand{\dual}{k[\epsilon]/\epsilon^{2}}
\newcommand{\ape}{A_{pe}}
\newcommand{\dtee} {
\otimes_{X}^{\mathbb{L}}
}
\newcommand{\dfm}[1]{\{\!\{\!#1\!\}\!\}}
\newcommand{\optp}{\mathscr{T}'^{op}}

\newcommand{\dteee} { \otimes^{\mathbb{L}} }
\newcommand{\ihom} {\mathbb{R}\underline{Hom}}
\newcommand{\A}{\mathscr{A}}

\newcommand{\F}{\mathscr{F}}
\newcommand{\G}{\mathscr{G}}
\newcommand{\K}{\mathscr{K}}

\newcommand{\M}{\mathscr{M}}
\newcommand{\Ox}{\mathscr{O}}
\newcommand{\Pp}{\mathscr{P}}
\newcommand{\T}{\mathscr{T}}
\newcommand{\R}{\mathscr{R}}
\newcommand{\Ss}{\mathscr{S}}

\newcommand{\der}[1][X]{D^{b}(#1)}
\newcommand{\opt}{\mathscr{T}^{op}}

\theoremstyle{plain}
\newtheorem{thm}{Theorem}[section]
\newtheorem*{thm*}{Theorem} 

\newtheorem{prop}[thm]{Proposition}
\newtheorem{lemma}[thm]{Lemma}
\newtheorem*{lemma*}{Lemma}
\newtheorem{defn}[thm]{Definition}
\newtheorem*{defn*}{Definition}

\newtheorem*{conjecture*}{Conjecture}
\newtheorem{exmp}{Example}[section]
\newtheorem{obs}[thm]{Remark}
\newtheorem{cor}[thm]{Corollary}
\newtheorem*{cor*}{Corollary}

\numberwithin{equation}{section}

\title{Davydov-Yetter cohomology for Tensor Triangulated Categories}
\author{Angel Toledo}
\email{toledo@unice.fr}
\date{}

\begin{document}
\begin{abstract}
    One way to understand the deformation theory of a tensor category $M$ is through its Davydov-Yetter cohomology $H_{DY}^{\ast}(M)$ which in degree 3 and 4 is known to control respectively first order deformations of the associativity coherence of $M$ and their obstructions. \\
    In this work we take the task of developing an analogous theory for the deformation theory of tensor triangulated categories with a focus on derived categories coming from algebraic geometry. We introduce the concept of perfect pseudo dg-tensor structure $\Gamma$ on an appropriate dg-category $\T$ as a truncated dg-lift of a tensor triangulated category structure on $\hzero(\T)$ and we define a double complex $DY^{\ast,\ast}(\Gamma)$ and we see that the 4th cohomology group $HDY^{4}(\Gamma)$ of the total complex of $DY^{\ast,\ast}(\Gamma)$ contains information about infinitesimal first order deformations of the tensor structure. 
\end{abstract}
\maketitle
\tableofcontents
\vspace{1in}

\include{Intro}
\include{Lifting}
\include{DYCohomology}

\hypersetup{linkcolor=red}
\bibliographystyle{alpha}
\bibliography{DDY.bib}

\end{document}

%% file: Intro.tex
\section{Introduction}
The study of tensor categories, nicely behaved k-linear abelian categories equipped with a rigid symmetric monoidal category structure, has been an active area of research in the past decades and has seen as a successful tool in representation theory ever since the work of Saavedra Rivano in Tannakian categories \cite{saavedra1972categories}. \\
The deformation theory of these categories was first introduced and studied by Davydov, Yetter and Crane-Yetter in \cite{davydov1997twisting,yetter1998braided,yetter2003abelian,crane1998deformations}. In concrete terms, this deformation theory is a deformation of the coherence data of the monoidal category or monoidal functor structures. Given a tensor category $(\A,\otimes)$, one defines a complex $DY^{\ast}(\A,\otimes)$ given by natural morphisms 
\[A_{1}\otimes (A_{2}\otimes (\dots \otimes(A_{n-1}\otimes A_{n})\dots)\to (\dots(A_{1}\otimes A_{2})\otimes \dots)\otimes A_{n-1})\otimes A_{n})\]
Where $A_{i}$ are possibly different objects of $\A$. \\
The cohomology of this complex is known to give information on the infinitesimal deformations of coherence conditions. In particular the 3rd cohomology group $HDY^{3}(\A,\otimes)$ controls the first order infinitesimal deformations of the associativity condition of the monoidal structure $\otimes$ and the group $HDY^{4}(\A,\otimes)$ controls the obstructions of this deformations (\cite[Exercise 7.22.2]{etingof2016tensor}). \\
If we denote by $\alpha_{A,B,C}:A\otimes(B\otimes C)\to (A\otimes B)\otimes C$ the natural associativity morphisms of the monoidal structure $\otimes$, by first order deformation of $\alpha$ we mean a natural family of morphisms \[
\alpha'_{A,B,C}:A\otimes'(B\otimes' C)\to (A\otimes' B)\otimes' C\]
Such that $\alpha'$ is the associativity data of a monoidal structure on a $\dual$-linear tensor category $\A'$ whose objects are the same as $\A$ and has $\otimes'$ as tensor bifunctor $\A'\times\A'\to \A'$, and such that $\alpha=\alpha'$ when restricted to $k$. \\
\\
In this work, we present a version of this constructions for tensor triangulated categories. To us, a tensor triangulated category is a $k$-linear triangulated category $\T$ equipped with a rigid symmetric monoidal category structure $\boxtimes$. In particular we should focus on tensor triangulated categories coming from algebraic geometry, meaning that the underlying triangulated category is the derived category of perfect complexes $Perf(X)$ on a nicely behaved scheme $X$. \\
Our motivation comes from an interest in tensor triangulated geometry as initiated by Balmer in \cite{balmer2002presheaves,balmer2005spectrum}, where a general construction is described to produce a locally ringed space $Spc(\T,\boxtimes)$ from a general tensor triangulated category $(\T,\boxtimes)$. While this space only depends on the unit object and the bifunctor in the monoidal structure on $\T$, we believe that there might be useful information on the categorical side of things which might reflect information which relates the geometry of $X$ with respect to its derived category of perfect complexes $Perf(X)$.\\
\\
In order to proceed with our constructions it is necessary to work with dg-enhancements of our triangulated categories. A dg-category is a category enriched over $\cc$. It is known that our derived categories $Perf(X)$ have unique dg-enhancements, meaning a dg-category $\T$ such that its homotopy category $\hzero(\T)$ is equivalent as a triangulated category to $Perf(X)$ where $\hzero(\T)$ is the category with the same objects as $\T$ and for any two objects $X,Y\in \hzero(\T)$, $Hom_{\hzero(\T)}(X,Y):=\hzero(Hom_{\T}(X,Y))$. \\
The reason we need to work at the level of dg-categories is that we will be exploiting Toën's Morita theorem from \cite{toen2007homotopy}. Recall that the classical Morita theorem completely characterizes additive cocontinuous functors between module categories $R-Mod$, $S-Mod$ as being equivalent to the category of $R-S$-bimodules. \\
We want then to characterize the exact bifunctors $\boxtimes:\T\times \T\to \T$ that form part of a tensor triangulated category structure $(\T,\boxtimes)$ in terms of certain bimodules. However, as it is well known, triangulated categories don't always behave sufficiently well and no reasonable version of this theorem holds for triangulated functors. \\
At the level of dg-categories however one such theorem is true when one takes into account equivalences of dg-categories up to weak equivalences in a certain model category structure on the category of dg-categories. \\
Inspired by \cite{hovey2011additive} where classical Morita theory was used to characterize additive closed symmetric monoidal category structures on a category of $R$-modules, we use Toën's Morita theorem to do the same for tensor triangulated category structures in terms of a structure we will call perfect pseudo dg-tensor structures on a given dg-category $\T$. \\
These structures can be thought of as a lift of tensor triangulated category structures to the dg-categorical level up to a given homotopy. \\
\\
We will proceed as follows: In section \ref{2} we will go over the general theory of dg-categories to properly estate the Morita theorem for dg-categories, then we will give a definition for pseudo dg-tensor category structures (\ref{def:pseudodgtensor}) and perfect pseudo dg-tensor category structures (\ref{def:perfectdgtensor}) by using n-fold dg-bimodules (\ref{def:nfoldbimodule}), we will then show that these structures induce tensor triangulated categories at the homotopy level (\ref{lemma:correspondence}) and give conditions for a bifunctor $\boxtimes$ part of a tensor triangulated category on a triangulated category $\T$ structure to lift to a perfect pseudo dg-tensor category structure on a dg-enhancement of $\T$. \\
Finally, in section \ref{3} we will define the Davydov-Yetter double complex associated to a perfect pseudo dg-tensor category structure $\Gamma$ (\ref{def:dgcohomologydavydovyetter}), the Davydov-Yetter cohomology of $\Gamma$, $HDY^{\ast}(\Gamma)$, and will prove our main theorem:
\begin{thm}(Theorem \ref{thm:deformationdg})
Let $\T$ be a dg-category and let $\Gamma$ be a pseudo dg-tensor structure on $\T$. Then to any element of  $HDY^{4}_{dg}(\T)$ we can associate an  equivalence class of infinitesimal deformations of order 1 of the associativity condition of $\Gamma$.
\end{thm}
Through this work unless stated otherwise, $k$ will denote a commutative ring, we denote by $\cc$ the category of cochain complexes on $k$. By scheme or variety we mean a scheme or variety over $k$.
\\
The results in this work were obtained as part of the author's PhD thesis at the Laboratoire J.A. Dieudonné at the Université Côte d'Azur. The author would like to thank his advisor Carlos Simpson for many discussions and to Ivo Dell'Ambrogio and Bertrand Toën for their careful and valuable comments on the thesis manuscript. The PhD thesis was partially financed by the CONACyT-Gobierno Francés 2018 doctoral scholarship.

%% file: Lifting.tex
\section{Lifting tensor triangulated categories}\label{2}
We will start this section by giving a quick review on the homotopy theory of dg-categories as developed by Toën in \cite{toen2007homotopy}. Our main goal is to be able to give a presentation on how the Morita theorem applies to derived categories of quasi-separated compact schemes. \\
Afterwards we will take inspiration from the work of Hovey in \cite{hovey2011additive} and present a characterization of a lift of a tensor triangulated category in term of bimodules at the dg-level by using Toën's Morita theorem.
Let us very briefly give the following series of definitions for a category of noncommutative derived schemes. This theory plays no substantial role in the following but it is an appropriate general category in which we can work. The reader can restrict themselves to working with dg-categories of perfect complexes coming from smooth projective varieties. 
\begin{defn}\label{defn:noncommuativescheme}
A (derived) noncommutative scheme $X$ over $k$ is a k-linear dg-category equivalent to $Perf_{dg}(R)$ where $R$ is a cohomologically bounded dg-algebra over $k$.
\end{defn}
A morphism of noncommutative schemes $X,Y$ is simply a quasi-functor $F:X\to Y$. Together with these morphisms we have a category of noncommutative schemes over $k$ which we denote by $NCSch_{k}$. \\
\begin{defn}
A noncommutative scheme $X=Perf_{dg}(R)$ is proper if $\bigoplus_{p\in \mathbb{Z}}H^{p}(Hom_{X}(M,N))$ are finite dimensional $k$-vector spaces for any two perfect modules $M,N\in X$.
\end{defn}
It can be shown that this is equivalent to the cohomology algebra $\bigoplus_{p\in \mathbb{Z}} H^{p}(R)$ being finite dimensional. \\
Furthermore, if $X$ is a separated scheme, then it is proper if and only if its dg-category $Perf_{dg}(X)$ is proper. \\
\begin{defn}
A noncommutative scheme $X$ is regular if its triangulated category has a strong generator.
\end{defn}
Similarly, when $X$ is a quasi-compact and separated scheme then it is regular as a noncommutative scheme if and only if it can be covered by open affine's $Spec(R_{i})$ with each $R_{i}$ being of finite global dimension. \\
\begin{defn}
A noncommutative scheme $X$ is smooth over $k$ if the bimodule $(x,y)\mapsto X(x,y)$ is a compact object of $X\otimes X^{op}$-Mod.
\end{defn}
As with the rest of the properties, in the case where $X$ is of finite type, a commutative scheme is smooth if and only if its category of perfect modules is smooth.  \\
Finally we have:
\begin{defn}
We say a dg-category is saturated if it is proper, smooth and pretriangulated and its triangulated category is idempotent complete. 
\end{defn}
A version of the story we have told so far about derived categories of smooth projective varieties regarding Serre functors, semi-orthogonal decompositions and even of Toën's Morita theorem can be reproduced for noncommutative schemes. We refer to \cite{orlov2018derived} for a more through summary of the theory developed so far. \\
\subsection{Preliminaries of dg-categories}
The contents of this subsection are merely expository, a more detailed outlook of the theory can be consulted for example in \cite{canonaco2017tour,keller2006differential,bertrand2011lectures}. \\
We start with some basic definitions:
\begin{defn}
A dg-category $\T$ is a category enriched over the closed monoidal category $\cc$, the category of cochain complexes over $k$. \\
Similarly, let $\T,\T'$ be dg-categories, a dg-functor $\F:\T\to \T'$ is a functor enriched over $\cc$.
\end{defn}
By considering enriched functors and enriched natural transformations we can then form a dg-category of functors $Fun_{dg}(\T,\T')$ between dg-categories $\T$ and $\T'$.
\begin{defn}
Let $\T$ be a dg-category. The homotopy category $\hzero(\T)$ of $\T$ is the category with the same class of objects as $\T$ and, if $x,y\in \hzero(\T)$ we let the Hom sets be given as 
\[ Hom_{\hzero(\T)}(x,y):= \hzero(Hom_{\T}(x,y)) \]
The homotopy functor $\hzero(\F):\hzero(\T)\to \hzero(\T')$ is the induced functor at the level of homotopy categories.
\end{defn}
Just as we do in categories, for any given dg-category $\T$ there exists an opposite dg-category $\T^{op}$ given by the same objects as those of $\T$ and whose Hom complexes are, for two objects $x,y\in \T^{op}$, the cochain complex $Hom_{\T}(y,x)$. \\
We are mainly interested in working with a module theory for dg-categories in the same way we do for regular or dg-algebras. \\
\begin{defn}\label{dgmoddef}
Let $\T$ be a dg-category, the category of (right) dg-modules over $\T$ is the dg-category of dg-functors $\T^{op}\to \cc$. We will denote this dg-category by $\T^{op}-Mod$. \\
Analogously we have a notion of left $\T$-module if we consider functors $\T\to \cc$. This category will be in turn denoted simply by $\T-Mod$.
\end{defn}
Let $\T$ be a dg-category and $x\in \T$ an object, then the module given by 
\[ h^{x}:=Hom_{\T}(x,\_):\T^{op}\to \cc \]
which associates $y\in \T$ the cochain complex $Hom_{\T}(x,y)$.
It is natural to consider the following definition
\begin{defn}\label{yonedahom}
Let $\T$ be a dg-category, the Yoneda embedding is the dg-functor:
\[ h^{\_}:\T\to \T^{op}-Mod \]
Which maps $x\in \T$ to $h^{x}$
\end{defn}
We can justify this nomenclature as it is possible to show that this dg-functor is full and faithful in the sense that there is an isomorphism of complexes for any $x,y\in \T$
\[ Hom_{\T^{op}-Mod}(h^{x},h^{y}) \cong Hom_{\T}(x,y) \].
We also have a co Yoneda embedding by the dg-functor
\[ h_{\_}:\T^{op}\to \T-Mod \]
Which takes $x\in \T$ and maps it to the left module $Hom_{\T}(\_,x):\T\to \cc$. \\
As usual we call a right $\T$-module representable if it is equivalent in $\T^{op}-Mod$ to a right module of the form $h^{x}$ for some $x\in T$. \\
As is the case for enriched categories, one can define a tensor product $\otimes$ of dg-categories which makes $dg-cat$ into a closed monoidal category. 
\begin{defn}\label{dgtensor}
If $\T,\T'$ are dg-categories, we define their tensor product as the dg-category $\T\otimes \T'$ with objects given by pairs $(x,x')\in \T\times \T'$, and cochain complexes between pairs of objects 
\[ Hom_{\T\otimes \T'}((x,x'),(y,y')):= Hom_{\T}(x,y)\otimes_{\cc} Hom_{\T'}(x',y') \]
And with composition here given entry-wise in the obvious way. \\
\end{defn}
Now that we are able to take these tensor products, we can consider a particular sort of dg-module which will play an essential role in our theory.
\begin{defn}
Let $\T, \T'$ be dg-categories. A left module  $\T\otimes \T'^{op}\to \cc$ will be called a dg-bimodule over $\T$ and $\T'$.
\end{defn}
Equivalences between enriched categories turn out to be too strong of a requirement for our goal of a Morita theorem for dg-categories. It is for this reason that there was a need to develop a homotopy theory for these categories so that we would obtain the right notion of equivalence. \\
This was done in the language of model category structures. In analogy to the simplicial category setting, we will now describe the Dwyer-Kan model category structure on the category of (small) dg-categories over $k$, $dg-cat_{k}$. \\
We will not only need to put a model structure on the whole category of dg-categories, but we will also be interested in putting model category structures on specific sorts of $\cc$-enriched categories in a way that makes the enrichment compatible with the model structure. It is in this interaction and compatibility between the global homotopy category of dg-categories and the internal model structure that we can put on a single specific $\cc$-enriched category where the core of some arguments will take place so it is important to describe both in some detail.
We are in particular interested in putting a model category structure on the dg-category of dg-modules over a dg-category $\T$. \\
This model structure on $\T$-Mod can be explicitly described by declaring a morphism $f:\F\to \G$ of $\T$-modules as a weak equivalence if for any $x\in \T$, the induced morphism $f_{x}:\F(x)\to \G(x)$ is a quasi-isomorphism. \\
We say that $f:\F\to \G$ is a fibration if the morphism $f_{x}$ is a fibration in $\cc$. \\
With this structure and using the fact that $\cc$ is a cofibrantly generated with generating cofibrations $I$ and generating trivial cofibrations $J$, then $\T$-Mod is also cofibrantly generated with generating cofibrations given by those morphisms of $\T$-modules of the form 
\begin{equation}\label{generatingcofibrations}
    h^{x}\otimes_{\cc} C \overset{id\otimes f}{\longrightarrow} h^{x}\otimes_{\cc} D
\end{equation}
Where $x,y\in \T$ and $f:C\to D\in I$. \\
As $\T$-Mod is a dg-category itself, the internal Hom object will simply be the Hom cochain complex objects of $\T$-Mod and the tensoring by $\cc$ is defined degree-wise. \\
Let us denote by $Ho(\M)$ the homotopy category of a model category $\M$, defined as the localization at the class of weak equivalences $W$ of $\M$, in other words $Ho(\M):=\M[W^{-1}]$. We should also consider the full category $Int(\M)$ of fibrant and cofibrant objects of $\M$. \\
One important thing to notice about dg-modules over a dg-category $\T$ is that the dg-modules $h^{x}\in \T$-Mod are fibrant as all cochain complexes $Hom_{\T}(y,x)$ are fibrant in the model category of $\cc$. \\
This implies that the Yoneda embedding of Definition \ref{yonedahom} factorizes through $Int(\T-\text{Mod}$). \\
\begin{defn}
Let $\T$ be a dg-category. A dg-module $\F\in \T^{op}$-Mod is called quasi-representable if it is equivalent in $Ho(\T^{op}-\text{Mod})$ to a module of the form $h_{x}$
\end{defn}
We denote the collection of such modules by $(\T-Mod)^{rqr}$. \\
We can now move on to describe the Dwyer-Kan model category structure on $dg-cat_{k}$.
\begin{defn}\label{dwyerkanmodel}
Let $\T, \T'$ be dg-categories over $k$. We say that a dg-functor $\F:\T\to \T'$ is a weak equivalence if
\begin{enumerate}
    \item It is quasi-fully faithful. This means that the induced cochain complex morphism $Hom_{\T}(x,y)\to Hom_{\T'}(\F(x),\F(y))$ is a quasi-isomorphism for all $x,y\in \T$.
    \item It is quasi-essentially surjective if $\hzero(\F):\hzero(\T)\to \hzero(\T')$ is essentially surjective. So for any $x'\in \hzero(\T')$, there exists $x\in \T$ such that $\F(x)\simeq x'$ in $\hzero(\T')$.
\end{enumerate}
A dg-functor $\F:\T\to \T'$ is on the other hand a fibration if
\begin{enumerate}
    \item The induced morphism of complexes $Hom_{\T}(x,y)\to Hom_{\T'}(\F(x),\F(y))$ is a fibration in the model structure of unbounded complexes. 
    \item For any isomorphism $u':x'\to y'\in \hzero(\T')$ and any $y\in \hzero(\T)$ such that $\F(y)=y'$, there is an isomorphism $u:x\to y$ in $\hzero(\T)$ such that $\hzero(\F)(u)=u'$.
\end{enumerate}
\end{defn}
These classes of morphisms in the category of dg-categories form a model category structure called the Dwyer-Kan model structure on $dg-cat_{k}$. The homotopy category $dg-cat_{k}[w^{-1}]$, where w is the class of weak equivalences above, of this structure will be denoted by $H_{qe}$. \\
We say that a dg-functor $f\in H_{qe}(\T,\T')$ between two dg-categories $\T,\T'$ is a quasi-functor. \\
Unfortunately the Dwyer-Kan model structure is not a monoidal model category with the tensor product we have discussed. It is however possible to construct a derived tensor product of dg-categories by using a cofibrant replacement. That is, we define a bifunctor \[\_\dteee\_:dg-Cat\times dg-Cat\to dg-Cat\]
defined by 
\[ \T\dteee \T' : = Q(\T)\otimes \T' \]
where $Q(\_)$ is a cofibrant replacement and $\otimes$ is the tensor product of dg-categories. \\
Importantly, the monoidal structure defined by this tensor product is closed and Toën gives the following characterization of the internal Hom object.
\begin{thm}
The monoidal category $(H_{qe}, \dteee)$ is closed and for any two dg-categories $\T,\T'$ there is a natural isomorphism in $H_{qe}$
\[ \ihom(\T,\T') \simeq Int(\T\otimes^{\mathbb{L}}\T'^{op}-Mod^{rqr}) \]
\end{thm}
In particular we have that
\[  \widehat{\T}:=Int(\T^{op}-Mod) \simeq \ihom(\T^{op},Int(\cc)) \]
in the homotopy category. \\
Alternatively, we might write $\T_{pe}$ instead of $\widehat{\T}_{pe}$ specially in the case when working over a field, where there is no need to consider the $(\widehat{-})$ operator, or when there is no chance of confusion. \\
Let us denote too by $\ihom_{c}(\widehat{\T},{\T'})$ the category of continuous dg-functors between the two dg-categories $\T, \T'$. By this we mean those functors $\widehat{\T}\to \T'$ such that they commute with direct sums when passing to the homotopy category $H_{qe}$. \\
\begin{defn}\label{def:perfectdgmodule}
We call $\ihom_{c}(\widehat{T},\widehat{T}')$ the dg-category of Morita morphisms from $\T\to \T'$, and perfect Morita morphisms the dg-category $\ihom(\widehat{T}_{pe},\widehat{T}'_{pe})$
\end{defn}
We can now formulate the main theorem of this section
\begin{thm}\cite[Theorem 7.1]{toen2007homotopy}\label{toenmorita}
Let $\T$ be a dg-category and let $y:\T\to \opt-Mod$ denote the Yoneda embedding. For any other dg-category $\Ss$, we have
\begin{enumerate}
    \item The pullback functor $y^{\ast}:\ihom_{c}(\widehat{\T},\widehat{\Ss})\to \ihom(\T,\widehat{\T})$ is an isomorphism in $H_{qe}$.
    \item The pullback functor $y^{\ast}:\ihom(\widehat{\T}_{pe},\widehat{\Ss}_{pe})\to \ihom(\T,\widehat{\Ss}_{pe})$ is an equivalence in $H_{qe}$
\end{enumerate}
\end{thm}
Using this result we arrive at the more well known form of the result, which can be thought as 
\begin{cor}
Let $\T$ and $\Ss$ be two dg-categories, then there exists a natural isomorphism in $H_{qe}$
\[ \ihom_{c}(\widehat{\T},\widehat{\Ss}) \simeq \widehat{\T^{op}\dteee \Ss} \] 
\end{cor}
Before we are able to continue with the concrete cases that interest us, let us quickly recall the notion of dg-enhancements for a triangulated category. 
\begin{defn}
Let $\K$ be a triangulated category. We say that a dg-category $\T$ is a dg-enhancement of $\K$ if there exists a triangulated equivalence 
\[ \epsilon:\hzero(\T)\to \K \]
\end{defn}
\begin{defn}
Let $\K$ be a triangulated category and let $\T$ and $\Ss$ two enhancements $\epsilon:\hzero(\T)\to \K, \epsilon':\hzero(\Ss)\to \K$. We say that $\K$ has a unique enhancement if there is a quasi-functor $f:\T\to\Ss$ such that $\hzero(f)$ is an equivalence of triangulated categories. 
\end{defn}
Importantly, let us consider $\T$ a dg-category and $\T$-Mod the dg-category of $\T$-modules. It can be shown then that the category $\hzero(\T-\text{Mod})$ can be equipped with a natural triangulated category structure. We then have:
\begin{defn}
Let $\T$ be a dg-category. We say that it is pretriangulated if its image under the functor $\hzero(y^{\T}):\hzero(\T)\to \hzero(\T-Mod)$ is a triangulated subcategory. 
\end{defn}
\begin{defn}\label{defn:pretriangulatedhull}
If $\T$ is a dg-category, we let $\T^{pre-tr}$ be the smallest pretriangulated full  dg-subcategory of $\T-Mod$. In this way we are adding cones, direct sums and all that might be missing from the original dg-category $\T$
\end{defn}. 
Let us denote the triangulated category $\hzero(\T^{pre-tr})$ by $tri(\T)$, and by $perf(\T)$ the full subcategory of compact objects in $tri(\T)$. \\
\begin{obs}
Given a dg-category $\T$, the triangulated category $perf(\T)$ is enhanced by $ \widehat{\T}_{pe}$.
\end{obs}
\begin{defn}
Let $\T$ be a dg-category and let $\Ss\subset \T$ be a full sub dg-category. The quotient $\T/\Ss$ is the dg-category with the same collection of objects as $\T$ and such that for every $s\in \Ss$ we add a morphism $s\to s$ in degree $-1$ so that $d(s\to s)=Id_{s}$
\end{defn}
With this definition of the quotient of dg-categories it is possible to show that $\hzero(\T)/\hzero(\Ss)\simeq \hzero(\T/\Ss)$, and we can now see that the derived category of an abelian category $\A$ can be enhanced by the quotient $C(\A)/Ac(\A)$.\\
The same can be done about bounded ( below, above, and both ) derived categories by simply taking the corresponding subcategories of $C(\A)$ and by consequence of $Ac(\A)$. \\
We can consider the derived category of a general dg-category $\T$ (cf. \cite{keller2006differential}). Let $Ac(\T)$ be the dg-subcategory of dg-modules $\T-Mod$ of those modules acyclic on every object. We have then:
\begin{defn}\label{defn:derivedcategorydgcategory}
Let $\T$ be a dg-category, the derived category $D(\T)$ is the quotient $\hzero(\T-Mod/Ac(\T))$.
\end{defn}
This derived category is always triangulated as it is a Verdier quotient of the triangulated category $\hzero(\T-Mod)$. \\
For a scheme we have that there always exists a dg-enhancement. The question of uniqueness can be deduced from Theorem B in \cite{canonaco2018uniqueness}.
Now that we have seen that our spaces and spaces of geometric nature can be enhanced in great generality, we are still left with the question of whether a triangulated functor can be lifted to a functor between dg-enhancements. \\
To be more explicit,
\begin{defn}
Let $\K, \K'$ be triangulated categories with dg-enhancements $\T, \T'$. We say that an exact functor $\F:\K\to \K'$ has a dg-lift if there exists a morphism $f\in H_{qe}(\T,\T')$ such that $\hzero(f)=F$.
\end{defn}
From now on let us fix and denote by $Perf_{dg}(X)$ a dg-enhancement of a derived category of perfect complexes over a space $X$, and similarly we put $QCoh_{dg}(X)$ for a dg-enhancement of the derived category of quasi coherent sheaves on $X. $\\
A consequence of Theorem \ref{toenmorita} is the following
\begin{thm}\label{thm:fouriermukaidg}
Let $X,Y$
 be two quasi-compact and separated schemes over k, assume that one is flat over $k$. Then there exists an isomorphism in $H_{qe}$
 \[ \ihom_{c}(QCoh_{dg}(X),QCoh_{dg}(Y))\simeq QCoh_{dg}(X\times_{k} Y) \]
 \end{thm}
 As a corollary of this, in the smooth case we have
\begin{thm}
Let $X$ and $Y$ be two smooth and proper schemes over k. Then there exists an isomorphism in $H_{qe}$
\[ \ihom(Perf_{dg}(X),Perf_{dg}(Y)) \simeq Perf_{dg}(X\times_{k} Y) \]
\end{thm}
Essentially what this theorem is telling us is that there is a correspondence in the smooth case between dg-lifts of exact functors between derived categories of perfect complexes and Fourier-Mukai transforms given by kernels in $Perf_{dg}(X\times_{k} Y)$. \\
Recall that a Fourier-Mukai transform $F:D^{b}(X)\to D^{b}(Y)$ is a functor equivalent to a functor of the form \[\mathbb{R}p_{\ast}({\mathbb{R}{q}}^{\ast}(\_\dteee K))\]
where $K\in D^{b}(X\times Y)$ and $p:X\times Y\to Y$, $q:X\times Y\to X$ are the projections to each factor. \\
\begin{obs}
In general, however, it is not true that exact functors between triangulated categories can be lifted not even in nice geometric settings. In \cite{rizzardo2019example} an explicit example is given of such exact functor that is not of Fourier-Mukai type. 
\end{obs}
Our focus in this work is entirely in this geometric setting rather than in arbitrary triangulated categories. It can be shown that the the dg-category $QCoh_{dg}(X)$ enhancing the derived category of quasi-coherent sheaves on a variety $X$ has a compact generator. This is an useful property as there is an equivalence in $H_{qe}$
\[ QCoh_{dg}(X)\simeq \widehat{A}_{X}\] 
where $A_{X}$ is a dg-algebra seen as a one-object dg-category.  \\
Using this characterization we can describe the object corresponding to 
\[ \F\in \ihom(Perf_{dg}(X),Perf_{dg}(Y)).\]
First, by \cite[Lemma 8.11]{toen2007homotopy} there is an equivalence in $H_{qe}$
\[ \ihom(A_{X},\widehat{A}_{Y}) \simeq \widehat{A_{X}\otimes_{k}^{\mathbb{L}} A_{Y}^{op}}\simeq QCoh_{dg}(X\times_{k} Y). \]
Then to any $\F$ as above, we associate the dg-bimodule $M\in \ihom(A_{X},\widehat{A}_{Y,pe})$ given by $M(\ast):=\F((\widehat{A}_{X,pe})(\ast,\_))$. \\
For example, if we let $\F:\ape\to \ape\in \ihom(\ape,\ape)$ where $A$ is a dg-algebra, the $\ape^{op}$-module given by $\F(A)$ has a left $A$-action since there is a chain complex morphism 
\[ A\to End(\F(A))\]
And as we can identify the representable module $\optp(y,\_)\simeq \F(A)$ with the object $y\in \optp$ which in turns has a right action by $A$, we obtain our right and left actions by $A$. \\
Furthermore we know that this bimodule induces a quasi-functor equivalent to $\F$ by the assignment \[M\mapsto \F(A)\otimes M\]
Where the tensor product is a tensor product of $\F(A)$ seen as a right module and $M\in A-Mod$ seen as left module. \\
This in turn means that there exists an isomorphism when passing to $\hzero(\ape)$ between $\F(A)\otimes M$ and $\F(M)$ for all $M\in \ape$. \\
With this idea in mind we will give a characterization of bimodules which produce tensor triangulated category structures at the homotopy level.
\subsection{Pseudo dg-tensor structures}
Let us take a moment to recall the general construction of tensor products between dg-modules. Let $\T,\R,\Ss$ be cofibrant ( in the Dwyer-Kan model structure ) dg-categories and let $\F\in \T\otimes\Ss^{op}-Mod$, $\G\in \Ss\otimes \R-Mod$. 

\begin{defn}\label{formula:productofdgmodules}We define the tensor product $\F\otimes_{\Ss} \G\in \T\otimes\R-Mod$ as, for any  $t\in \T$, $r\in \R$, the chain complex calculated as the cokernel of 
\begin{equation*} \bigoplus_{x,y\in\Ss}  \F(t,x)\otimes_{k} \Ss(y,x) \otimes_{k} \G(y,r) \to \bigoplus_{z\in \Ss} \F(t,z)\otimes_{k} \G(z,r),\end{equation*}
where the morphism takes a homogeneous element $v\in \F(t,x)$, an homogeneous element $u\in \G(y,r)$ and a homogeneous morphism $f\in \Ss(y,x)$ to the homogeneous element 
\[ \F(t,f)(v)\otimes u - (-1)^{\mid v\mid \mid u \mid} v\otimes\G(f,r)(u) \]
With this construction in mind we define the following notion of n-fold dg-bimodules over a dg-category $\T$ 
\end{defn}
\begin{defn}\label{def:nfoldbimodule}
Let $\T$ be a dg-category. An n-fold dg-bimodule over $\T$ is a dg-module $\F\in \T^{\otimes n}\otimes \T^{op}-Mod$.
\end{defn}
In particular a 0-fold dg-bimodule is nothing but a $\T^{op}-\text{module}$ and a 1-fold bimodule is what we usually call a bimodule over $\T$. \\
A morphism of n-fold dg-bimodules is simply a morphism of dg-modules and we have then a dg-category denoted by $Bimod^{n}_{dg}(\T)$ with n-fold dg-bimodules as objects and morphism objects given by morphisms of dg-modules.  \\
Notice that the permutation group $\Sigma_{n}$ acts on $Bimod^{n}_{dg}$ by switching the tensor product $\T\;^{\otimes n}$ and so if $\F\in Bimod^{n}_{dg}(\T)$, $\sigma \in \Sigma_{n}$ is a permutation and $x\to y\in \T$ where this $\T$ lies in the k-th slot of the product $\T\otimes\dots\otimes\T$, $x\to y$ induces a morphism of dg-modules 
\begin{equation*}
    \begin{xymatrix}{ \F(\_,\dots,\underbrace{x}_{k},\dots,\_,?)\ar[r] \ar[d] & \sigma\F(\_,\dots,\underbrace{x}_{\sigma(k)},\dots,\_,?) \ar[d] \\ \F(\_,\dots, \underbrace{y}_{k},\dots,\_,?)\ar[r] & \sigma\F(\_,\dots,\underbrace{y}_{\sigma(k)},\dots,\_,?)
    }
    \end{xymatrix}
\end{equation*}
So that $x,y$ now lie in the $\sigma(k)$-th slot of $\sigma\F$. \\
To avoid confusion and ease the reading we follow Hovey's notation and introduce dummy variables to keep track of which $\T$ factor is being taken into account. For example if $\F$ and $\G$ are 3-fold dg-bimodules and $\sigma=(3 1)$, a morphism $\eta:\F\to \sigma\G$ is better expressed as 
\[ \eta:\F_{X,Y,Z}\to \G_{Z,Y,X} \]
to indicate that the action of any morphisms $x\to y\in \T$ at the first slot is now carried to an action to the third one on $\G$. \\
Using the tensor product of bimodules described in \ref{formula:productofdgmodules} we see that there exists, for any pair of natural numbers $n,m$ a way to tensor n-fold bimodules with m-fold bimodules. \\
If $\F\in Bimod_{dg}^{n}(\T)$ and $\G\in Bimod_{dg}^{m}(\T)$ then we form the tensor product $\F\otimes_{\T}\G$ by using the leftmost $\T$ factor in $\G$ with the right $\T$ factor of $\G$. \\
In the case we would want to take this tensor product with any of the other $\T$ right factors of $\F$ we can simply consider a permutation $\sigma$ which permutes the n-th factor with the factor we want to tensor with.\\
When using the notation above we are able to drop $\sigma$ from our expression as it is implied from the order of the subindices which permutation we are operating.  \\
To keep track of which factor is being used to form the tensor product we will extend the notation for morphisms and write for a 2-fold dg-bimodule $\F$, a 3-fold dg-bimodule $\G$
\[\F_{\G,X}\otimes \G_{Y,Z,W} \]
To denote we are forming a 4-fold dg-bimodule by using the first left factor of $\G$ to produce the tensor product with the right factor of $\F$.
In our particular case we are for the moment only interested in categories of modules over the a dg-algebra $A$ seen as a single object dg-category and so what we are describing is simply the theory of dg-bimodules over a dg-algebra with multiple compatible left multiplications and the tensor product described here is just the usual tensor product of right and left modules.  \\
We will make use of  Theorem \ref{toenmorita} to show that if have a functor $\F:\K^{n}\to \K$ from a product of a triangulated category $\K$ such that it is a triangulated functor in each variable, then it is possible to find an appropriate quasi-functor and make it correspond an n-fold dg-bimodule. \\
Having established our notation we now present a homotopical dg-module version of Theorem 2.1 in \cite{hovey2011additive}
\begin{thm}\label{thm:ttcasbimod}
Let $A$ be a dg-algebra and let $\boxtimes:\hzero(\ape)\times\hzero(\ape)\to \hzero(\ape)$ be an exact functor in each variable. Suppose that for every object $M\in \hzero(\ape)$, the triangulated functors
\begin{align*}
    M\boxtimes \_ :\hzero(\ape) \to \hzero(\ape)
    \\
    \_\boxtimes M:\hzero(\ape) \to \hzero(\ape)
\end{align*}
both have unique dg-enhancements $R_{M}$ and $L_{M}$ respectively. \\
Then $L_{A}(A)$ is a 2-fold dg-bimodule and for any $N\in \ape$ we have
\[ \hzero(L_{A}(A)\otimes M\otimes N)\simeq M\boxtimes N \]

\begin{proof}
We have to make repeated use of Theorem \ref{toenmorita}. First let $R_{M}\in \ihom(\ape,\ape)$ be an enhancement of $M\boxtimes \_$, we know by the derived Morita theorem that this quasi-functor corresponds to a bimodule given by $R_{M}(A)$. We know by the theorem that for any $N$, $\hzero(R_{M}(A)\otimes N)\simeq M\boxtimes N$. \\
Now consider the functor $\_\boxtimes A$ which by hypothesis has a unique enhancement $L_{A}$, by using Theorem \ref{toenmorita} again, we know this quasi-functor corresponds to a dg-bimodule $L_{A}(A)$ and that for any M we have that $\hzero(L_{A}(A)\otimes M)\cong M\boxtimes A$. \\
However, this latter object is isomorphic to $\hzero(R_{M}(A)\otimes A)$ and as this is a natural isomorphism on $M$ then the functor $R_{\_}(A)\otimes A$ defined by $M\mapsto R_{M}A\otimes A$ is too an enhancement of $\_\boxtimes A$. By uniqueness of the enhancement of $\_\boxtimes A$ we have then $L_{A}(A)\otimes \_ \simeq R_{\_}(A)\otimes A$ which corresponds to a right quasi-representable $A$-bimodule. \\
However, $R_{M}(A)$ is already a bimodule for every $M$, in other words $R_{A}(A)$ has two compatible dg-bimodule structures (in the sense these two actions on the same side by $A$ commute) and so we can consider as an object in $Bimod_{dg}^{2}(A_{pe})$.
Finally we obtain the required isomorphism:
\[\hzero(L_{A}(A)\otimes M\otimes N) \cong \hzero(R_{M}(A)\otimes N) \cong M\boxtimes N \]
\end{proof}
\end{thm}
Even though in general the existence of dg-lifts is not guaranteed as we mentioned before, in our particular context the existence of the lift is only a mild assumption. Indeed as we are coming from a dg-algebra induced from perfect complexes over a smooth projective variety, in that case it is  a consequence of the derived Morita theorem that a triangulated functor is of Fourier-Mukai type if and only if it has a dg-lift. \\
While the condition is a necessary one, as it is unknown whether every triangulated functor that we consider is of Fourier-Mukai type, we find this condition not too strong to impose. \\
The uniqueness of the lift, however is a stronger condition as we know for sure there are examples of non-uniqueness of the integral kernels that determine these transforms. \\
With this result in mind, we are finally off to approach tensor triangulated categories as being induced by 2-fold bimodules together with structure maps and their corresponding coherence conditions. Before going further let us remark the following
\begin{obs}
Since any exact functor between derived categories of smooth projective varieties has a right and left adjoints, if $\boxtimes$ is a tensor product of a tensor triangulated category in such triangulated category, the hypothesis of exactness in each variable implies automatically that the symmetric monoidal structure is closed. Let us remark that any of the triangulated functors that show up in our context are $k$-linear. \\
\end{obs}
Let us present the following definition in obvious analogy with the usual (lax) symmetric monoidal category axioms:
\begin{defn}\label{def:pseudodgtensor}
A pseudo dg-tensor structure in a dg-category $\T$ consists on the data:
\begin{enumerate}
    \item A 2-fold dg-bimodule $\Gamma \in Bimod^{2}_{dg}(\T)$
    \item An object $U\in \opt-Mod$ called the unit.
    \item  Morphisms of dg-bimodules $\alpha_{X,Y,Z}:\Gamma_{X,\Gamma}\otimes \Gamma_{Y,Z}\to \Gamma_{\Gamma,Z}\otimes \Gamma_{X,Y}\in Bimod_{dg}^{3}(\T)$ .
    \item A  morphism of dg-bimodules $\ell_{X}:\Gamma_{U,X}\otimes U \to \T-Mod \in Bimod^{1}_{dg}(\T)$.
    \item A  morphism of dg-bimodules $r_{X}:\Gamma_{X,U}\otimes U\to \T-Mod \in Bimod^{1}_{dg}(\T)$.
    \item A morphism $c_{X,Y}:\Gamma_{X,Y} \to \Gamma_{Y,X}$ of dg-bimodules.
\end{enumerate}
We require that the morphisms $\alpha_{X,Y,Z}$, $u_{X}$ and $c_{X,Y}$ are all isomorphisms when passing to the homotopy category $\hzero(\cc)$, for all $X,Y,Z\in \opt-Mod$.  \\
Furthermore we require the following homotopy data satisfying conditions:
\begin{enumerate}
    \item  (Associativity) A morphism $\eta\in Hom^{-1}(\Gamma_{X,\Gamma}\otimes \Gamma_{Y,\Gamma}\otimes \Gamma_{Z,W}, \Gamma_{\Gamma,W}\otimes \Gamma_{\Gamma,Z}\otimes \Gamma_{X,Y} )$ such that  $\alpha^{0}_{\Gamma_{X,Y},Z,W}\circ \alpha^{0}_{X,Y,\Gamma_{Z,W}}-\alpha^{0}_{X,Y,Z}\otimes Id_{W}\circ \alpha^{0}_{X,\Gamma_{Y,Z},W}\circ Id_{X}\otimes \alpha_{Y,Z,W}=d(\eta)$.
    \item (Unit) A morphism $\mu\in Hom^{-1}(\Gamma_{X,\Gamma}\otimes \Gamma_{U,Y}, \Gamma_{\Gamma,Y}\otimes \Gamma_{X,U}) $ such that  $\ell^{0}_{X}\otimes Id_{Y}\circ \alpha^{0}_{X,U,Y}-Id_{X}\otimes \ell^{0}_{Y}=d(\mu)$
    \item (Symmetry) The composition $c_{X,Y}\circ c_{Y,X}$ is the identity in $\hzero(\T-Mod)$. 
    \item (Unit symmetry) There is $\kappa\in Hom^{-1}(\Gamma_{X,U}, X)$ such that $\ell_{X}\circ c_{X,U}-r_{X}=d(\kappa)$
    \item (Compatibility between associativity and symmetry) There is $\lambda\in Hom^{-1}(\Gamma_{\Gamma,Z}\otimes \Gamma_{X,Y}, \Gamma_{Y,\Gamma}\otimes \Gamma_{Z,X})$ such that $Id_{Y}\otimes c_{X,Z}\circ \alpha_{X,Y,Z} \circ c_{X,Y}\otimes Id_{Z} - \alpha_{Y,Z,X} \circ c_{X,\Gamma_{Y,Z}} \circ \alpha_{X,Y,Z} = d(\lambda)$.
\end{enumerate}
\end{defn}
Here $\T-Mod$ denotes the denotes the $\T-\T$-bimodule 
\[ (x,y)\mapsto Hom_{\T-Mod}(y,x). \]
Our requirement that structure morphisms $\alpha,u, c$ are invertible in the homotopy category implicitly implies that the degree 0 part of these morphisms is a cycle. \\
These are structures up to 
As the coherence conditions are troublesome to keep track, let us use our graphical notation to draw the diagrams we require \\
\begin{center}\label{figure:associativitysymmetry}
\scalebox{0.5}{
\begin{tangle}{(18,10)}
	\tgBorderC{(2,0)}{3}{white}{white}
	\tgBlank{(3,0)}{white}
	\tgBorderC{(1,1)}{3}{white}{white}
	\tgBorderA{(2,1)}{white}{white}{white}{white}
	\tgBorder{(2,1)}{1}{0}{1}{1}
	\tgBorderC{(8,1)}{3}{white}{white}
	\tgBlank{(9,1)}{white}
	\tgBorderC{(14,1)}{3}{white}{white}
	\tgBlank{(15,1)}{white}
	\tgBorderC{(16,1)}{3}{white}{white}
	\tgBlank{(17,1)}{white}
	\tgBorderA{(0,2)}{white}{white}{white}{white}
	\tgBorder{(0,2)}{0}{1}{0}{0}
	\tgBorderA{(1,2)}{white}{white}{white}{white}
	\tgBorder{(1,2)}{1}{0}{1}{1}
	\tgBorderC{(2,2)}{0}{white}{white}
	\tgBorderC{(3,2)}{3}{white}{white}
	\tgBlank{(4,2)}{white}
	\tgBorderA{(5,2)}{white}{white}{white}{white}
	\tgBorder{(5,2)}{0}{1}{0}{1}
	\tgBorderA{(6,2)}{white}{white}{white}{white}
	\tgBorder{(6,2)}{0}{1}{0}{1}
	\tgBorderA{(7,2)}{white}{white}{white}{white}
	\tgBorder{(7,2)}{0}{1}{0}{0}
	\tgBorderA{(8,2)}{white}{white}{white}{white}
	\tgBorder{(8,2)}{1}{0}{1}{1}
	\tgBorderC{(9,2)}{3}{white}{white}
	\tgBlank{(10,2)}{white}
	\tgBorderA{(11,2)}{white}{white}{white}{white}
	\tgBorder{(11,2)}{0}{1}{0}{1}
	\tgBorderA{(12,2)}{white}{white}{white}{white}
	\tgBorder{(12,2)}{0}{1}{0}{1}
	\tgBorderA{(13,2)}{white}{white}{white}{white}
	\tgBorder{(13,2)}{0}{1}{0}{0}
	\tgBorderA{(14,2)}{white}{white}{white}{white}
	\tgBorder{(14,2)}{1}{0}{1}{1}
	\tgBorderC{(15,2)}{3}{white}{white}
	\tgBorderA{(16,2)}{white}{white}{white}{white}
	\tgBorder{(16,2)}{1}{0}{1}{1}
	\tgBorderC{(1,3)}{0}{white}{white}
	\tgBorderA{(2,3)}{white}{white}{white}{white}
	\tgBorder{(2,3)}{0}{1}{0}{1}
	\tgBorderA{(3,3)}{white}{white}{white}{white}
	\tgBorder{(3,3)}{1}{0}{1}{1}
	\tgBorderC{(8,3)}{0}{white}{white}
	\tgBorderA{(9,3)}{white}{white}{white}{white}
	\tgBorder{(9,3)}{1}{0}{1}{1}
	\tgBorderC{(10,3)}{3}{white}{white}
	\tgBlank{(11,3)}{white}
	\tgBorderC{(14,3)}{0}{white}{white}
	\tgBorderA{(15,3)}{white}{white}{white}{white}
	\tgBorder{(15,3)}{1}{0}{1}{1}
	\tgBorderC{(16,3)}{0}{white}{white}
	\tgBlank{(17,3)}{white}
	\tgBorderA{(2,4)}{white}{white}{white}{white}
	\tgBorder{(2,4)}{1}{0}{1}{0}
	\tgBorderC{(3,4)}{0}{white}{white}
	\tgBlank{(4,4)}{white}
	\tgBorderC{(9,4)}{0}{white}{white}
	\tgBorderA{(10,4)}{white}{white}{white}{white}
	\tgBorder{(10,4)}{1}{0}{1}{1}
	\tgBorderC{(15,4)}{0}{white}{white}
	\tgBorderA{(2,5)}{white}{white}{white}{white}
	\tgBorder{(2,5)}{1}{0}{1}{0}
	\tgBorderC{(3,5)}{3}{white}{white}
	\tgBlank{(4,5)}{white}
	\tgBorderC{(10,5)}{0}{white}{white}
	\tgBlank{(11,5)}{white}
	\tgBorderA{(15,5)}{white}{white}{white}{white}
	\tgBorder{(15,5)}{1}{0}{1}{0}
	\tgBorderC{(2,6)}{3}{white}{white}
	\tgBorderA{(3,6)}{white}{white}{white}{white}
	\tgBorder{(3,6)}{1}{0}{1}{1}
	\tgBlank{(14,6)}{white}
	\tgBorderC{(15,6)}{3}{white}{white}
	\tgBlank{(16,6)}{white}
	\tgBorderC{(1,7)}{3}{white}{white}
	\tgBorderA{(2,7)}{white}{white}{white}{white}
	\tgBorder{(2,7)}{1}{0}{1}{1}
	\tgBorderC{(3,7)}{0}{white}{white}
	\tgBlank{(4,7)}{white}
	\tgBorderA{(5,7)}{white}{white}{white}{white}
	\tgBorder{(5,7)}{0}{1}{0}{1}
	\tgBorderA{(6,7)}{white}{white}{white}{white}
	\tgBorder{(6,7)}{0}{1}{0}{1}
	\tgBorderA{(7,7)}{white}{white}{white}{white}
	\tgBorder{(7,7)}{0}{1}{0}{1}
	\tgBorderA{(8,7)}{white}{white}{white}{white}
	\tgBorder{(8,7)}{0}{1}{0}{1}
	\tgBorderA{(9,7)}{white}{white}{white}{white}
	\tgBorder{(9,7)}{0}{1}{0}{1}
	\tgBorderA{(10,7)}{white}{white}{white}{white}
	\tgBorder{(10,7)}{0}{1}{0}{1}
	\tgBorderA{(11,7)}{white}{white}{white}{white}
	\tgBorder{(11,7)}{0}{1}{0}{1}
	\tgBorderA{(12,7)}{white}{white}{white}{white}
	\tgBorder{(12,7)}{0}{1}{0}{1}
	\tgBlank{(13,7)}{white}
	\tgBorderC{(14,7)}{3}{white}{white}
	\tgBorderA{(15,7)}{white}{white}{white}{white}
	\tgBorder{(15,7)}{1}{0}{1}{1}
	\tgBorderC{(16,7)}{3}{white}{white}
	\tgBlank{(17,7)}{white}
	\tgBorderA{(0,8)}{white}{white}{white}{white}
	\tgBorder{(0,8)}{0}{1}{0}{0}
	\tgBorderA{(1,8)}{white}{white}{white}{white}
	\tgBorder{(1,8)}{1}{0}{1}{1}
	\tgBorderC{(2,8)}{0}{white}{white}
	\tgBlank{(3,8)}{white}
	\tgBorderA{(13,8)}{white}{white}{white}{white}
	\tgBorder{(13,8)}{0}{1}{0}{0}
	\tgBorderA{(14,8)}{white}{white}{white}{white}
	\tgBorder{(14,8)}{1}{0}{1}{1}
	\tgBorderC{(15,8)}{0}{white}{white}
	\tgBorderA{(16,8)}{white}{white}{white}{white}
	\tgBorder{(16,8)}{1}{0}{1}{1}
	\tgBorderC{(1,9)}{0}{white}{white}
	\tgBlank{(2,9)}{white}
	\tgBorderC{(14,9)}{0}{white}{white}
	\tgBorderC{(16,9)}{0}{white}{white}
	\tgBlank{(17,9)}{white}
	\tgCell{(8.5,1)}{X}
	\tgCell{(9.5,2)}{Y}
	\tgCell{(10.5,3)}{Z}
	\tgCell{(10.5,5)}{W}
	\tgCell{(3.5,2)}{Z}
	\tgCell{(2.5,2)}{Y}
	\tgCell{(3.5,4)}{W}
	\tgCell{(3.5,7)}{Y}
	\tgCell{(3.5,5)}{X}
	\tgCell{(2.5,8)}{Z}
	\tgCell{(1.5,9)}{W}
	\tgCell{(2.5,0)}{X}
	\tgArrow{(12.5,2)}{0}
	\tgArrow{(2,5.5)}{3}
	\tgArrow{(4.5,7)}{2}
	\tgArrow{(15,5.5)}{3}
	\tgCell{(15.5,6)}{X}
	\tgCell{(16.5,7)}{Y}
	\tgCell{(15.5,4)}{W}
	\tgCell{(14.5,1)}{X}
	\tgCell{(16.5,1)}{Y}
	\tgCell{(16.5,3)}{Z}
	\tgCell{(14.5,9)}{W}
	\tgCell{(16.5,9)}{Z}
	\tgArrow{(4.5,2)}{2}
	\tgCell[(1.5,0)]{(5.7,2)}{\alpha_{X,Y,\Gamma_{Z,W}}}
	\tgCell[(1.9,0)]{(11.2,2)}{id\otimes \alpha_{Y,Z,W}}
	\tgCell[(2,0)]{(8.5,7)}{\alpha_{X,Y,Z}\otimes id_{W}}
	\tgCell[(1.5,0)]{(2,4.5)}{\alpha_{\Gamma_{X,Y},Z,W}}
	\tgCell[(1.5,0)]{(15,5)}{\alpha_{X,\Gamma_{Y,Z},W}}
        \tgCell[(1,0)]{(8,5)}{d(\eta)}
\end{tangle}
}
\end{center}
This is nothing but the pentagon axiom in a monoidal category except we don't require this composition to commute but only to commute up to the homotopy $\eta$, $d(\eta)$. \\
Similarly, we can obtain analogous diagrams corresponding to the rest of the coherence conditions, we hope that the associativity diagram above sufficiently explain the equations from Definition \ref{def:pseudodgtensor}.\\
We should interpret this diagram as each branch in the tree coming from tensoring by $\Gamma$ and letting us tensor by a single left $\T$-module on each end. \\
\begin{obs}
We could alternatively present the previous construction from an operadic point of view. See for example \cite[Chapter 1]{beilinson2004chiral} for a presentation of close ideas in this language. The reader should keep in mind that an incompatibility in nomenclature. 
\end{obs}
We will show that given a pseudo dg-tensor structure on a dg-category $\T$, we can induce a tensor triangulated category structure on $\widehat{T}_{pe}$ under certain conditions. \\
By abuse of notation and when there is no chance for confusion, we denote by $\Gamma$ the pseudo dg-tensor structure ($\Gamma,U,\alpha,\ell, r,c)$. \\
An immediate question to consider is how to declare two such structures as equivalent. Let us say that two such pseudo dg-tensor category structures $\Gamma,\Gamma'$ are \textit{tensor triangulated equivalent} if they induce the same tensor triangulated structure on $\hzero(\T-Mod)$. This is a somewhat weak notion as it is entirely possible for two non equivalent 2-fold dg-bimodules to produce this same tensor triangulated structures when passing to the homotopy category, or even just the defining morphisms can be different at higher degrees while inducing the same structure homotopically. This is the reason this lift to the dg-setting has to be thought of as a truncated version.  \\
Let us be more precise in our claims above, first we need to introduce the following refinement of our definition above: \\
\begin{defn}\label{def:perfectdgtensor}
A pseudo dg-tensor structure $\Gamma$ in a dg-category $\T$ is called perfect if the 2-fold dg-bimodule $\Gamma\in Bimod_{dg}^{2}(\T)$ is right quasi-representable and for every $X,Y\in \opt-Mod$, $\Gamma\otimes X\otimes Y$ is quasi-represented by a perfect $\opt$-module.
\end{defn}
\begin{lemma}\label{lemma:correspondence}
A perfect pseudo dg-tensor structure $\Gamma$ on a dg-category $\T$ induces a tensor triangulated category structure on $\hzero(\widehat{\T}_{pe})$.
\begin{proof}
The functor $\hzero(\Gamma_{X,Y}\otimes\_\otimes\_)$ induces a bifunctor on $\hzero(\T_{pe})$. Indeed let us denote by $X\boxtimes Y$ the equivalence class of perfect modules quasi-representing $\Gamma_{X,Y}\otimes X\otimes Y$, by hypothesis we know $\hzero(X\boxtimes Y)$ is a perfect object and fixing $X$ or $Y$ we get a quasi-representable dg-bimodule which induces triangulated functors
\[ \_\boxtimes Y:\hzero(\T_{pe})\to \hzero(\T_{pe}) \]
and
\[ X\boxtimes \_:\hzero(\T_{pe})\to \hzero(\T_{pe}) \]
Using the dg-bimodule morphisms $\alpha_{X,Y,Z}$, $\ell_{X}$, $r_{X}$, $c_{X,Y}$ we obtain morphisms $X\boxtimes (Y\boxtimes Z)\to (X\boxtimes Y)\boxtimes Z$, $U\boxtimes X \to X$, $X\boxtimes U\to X$, and $X\boxtimes Y \to Y\boxtimes X$. \\
The equations in Definition \ref{def:perfectdgmodule} encode the coherence conditions for this monoidal category with product functor $\boxtimes$. For example, the associativy can be seen from the pentagon diagram \ref{figure:associativitysymmetry} as explained above. Indeed, the condition that they must commute up homotopy means that when passing to the homotopy category $\hzero(\widehat{T}_{pe})$ they will commute in the usual sense. 
\end{proof}
\end{lemma}
We have to remark a couple things. The first one is that for us, it is necessary to include both right and left unit conditions in the pseudo dg-tensor structure even as if the existence of $c_{X,Y}$ implies that one can obtain one from the other in the 1-categorical setting. In our case however it is necessary to keep track of them as separate entities. \\
The second thing to mention is that seeing our lemma as a dg-version of Theorem 2.3 of \cite{hovey2011additive} we need to remark that the converse does not hold as-is. Indeed as liftings of objects and morphisms in a triangulated category to a dg-enhancement are far from being unique we cannot expect to have unique -up to isomorphism- dg-bimodules, structure maps and homotopies inducing a certain tensor triangulated category. 
As in the abelian case we can too encode what a lift of a triangulated tensor endofunctor would be in our dg-setting, and as such a morphism between perfect pseudo dg-tensor structures on a given dg-category $\T$.
\begin{defn}\label{def:pseudodgtensorfunctor}
Let $\T$ be a dg-category and let $\Gamma$ and $\Lambda$ be perfect pseudo dg-tensor structures on $\T_{pe}$ with units $U$ and $U'$ respectively. A pseudo dg-tensor functor between $\Gamma$ and $\Lambda$ consists of,
\begin{enumerate}
    \item A dg-bimodule $\Phi$ 
    \item A morphism of dg-modules $u:U\to \Phi\otimes U'$
    \item A morphism of dg-bimodules $f:\Gamma_{\Phi, \Phi}\otimes\Phi\otimes\Phi \to \Phi\otimes \Lambda$ 
\end{enumerate}
Such that these two morphisms are isomorphisms when passing to $\hzero$. \\
Furthermore, we need the following coherence conditions
\begin{enumerate}
    \item There exists $w\in Hom^{-1}(\Gamma_{\Phi,U'}\otimes \Phi\otimes U, X)$ such that $Id_{\Phi}\otimes \ell \circ f\circ Id_{\Gamma}\otimes Id_{\Phi}\otimes T_{U,\Phi}\circ Id_{\Gamma}\otimes u \otimes Id_{\Phi}- \ell\otimes Id_{X}=d(w)$
    \item There exists $e\in Hom^{-1}(\Gamma_{\Phi,\Phi}\otimes \Phi\otimes \Phi, \Phi\otimes \Lambda)$ such that $ Id_{\Phi}\otimes c_{\Lambda}\circ f- f\circ c_{\Gamma}\otimes T_{\Phi,\Phi}=d(e)$
    \item There exists $a\in Hom^{-1}(\Gamma_{\Phi,\Gamma}\otimes \Gamma_{\Phi,\Phi}\otimes \Phi\otimes\Phi\otimes\Phi, \Phi\otimes\Lambda\otimes\Lambda)$ such that $f\circ T_{\Phi\otimes \Lambda, \Phi} \circ Id_{\Gamma}\otimes f\otimes Id_{\Phi}\circ \alpha_{\Phi,\Phi,\Phi}\otimes T_{\Phi,\Phi}-f\circ T_{\Phi\otimes \Lambda, \Phi}\circ f\otimes Id_{\Lambda}\otimes Id_{\Phi}\otimes \alpha= d(a)$
\end{enumerate}
\end{defn}
The structure morphisms and coherence conditions above in the definition are nothing but the structural morphisms and coherence conditions of a monoidal functor with the equivalent underlying category written in terms of bimodules and morphisms between them, with the only difference as in the structure of a pseudo dg-tensor structure being that we have to specify a given homotopy. \\
The proof of the following is straightforward:
\begin{lemma}\label{lemma:dgtensorfunctorbimod}
Let $\T$ be a dg-category and let $\Gamma$ and $\Lambda$ be  perfect pseudo dg-tensor structures on $\T_{pe}$. Then a pseudo dg-tensor functor $\Phi$ from $\Gamma$ to $\Lambda$ induces a tensor triangulated functor $\F_{\Phi}:\hzero(\T_{pe},\Gamma)\to\hzero(\T_{pe},\Lambda)$, where $\hzero(\T_{pe},\Gamma)$ and $\hzero(\T_{pe},\Lambda)$ denote the tensor triangulated categories induced by $\Gamma$ and $\Lambda$ respectively 
\begin{proof}
We saw from Lemma \ref{lemma:correspondence} that $\Gamma$ and $\Lambda$ produces tensor triangulated structures on $\hzero(\T_{pe})$ and by Toën's Morita theorem, a triangulated functor $\T_{pe}\to \T_{pe}$ corresponds to a $\T$-bimodule $\Phi$.  The structural morphisms $u$ and $f$ of $\Phi$ being isomorphisms in $\hzero(\T_{pe})$, and the coherence conditions imply that the induced functor between triangulated categories is a symmetric monoidal functor.
\end{proof}
\end{lemma}
Composition of functors corresponds to tensor product of bimodules. There exists a canonical identity pseudo dg-tensor functor which is given by the dg bimodule $\T-Mod$ (meaning the dg-bimodule which sends $X\in \T^{op}, Y\in\T$ to $\T(X,Y)$) and structural and coherence conditions all given by the canonical isomorphisms $U\to \T-Mod\otimes U'$ and $\Gamma\otimes \T-Mod\otimes\T-Mod \to \T-Mod\otimes\Gamma$.\\
Two pseudo dg-tensor functors $\Phi$ and $\Phi'$ are said to be equivalent if there exists a morphism of bimodules $\Phi\to \Phi'$ such that $\hzero(\Phi)\to \hzero(\Phi')$ is an isomorphism and is compatible with the morphisms $u$ and $f$ in the obvious homotopical sense. \\
We do not describe these natural transformations in detail as we will not be needing the coherence conditions of them but only knowledge that the usual diagrams commute up to a given homotopy. \\
Previously we said that two perfect pseudo dg-tensor structures were tensor triangulated equivalent if they gave rise to equivalent tensor triangulated structures.  We have the following lemma that compares the two equivalence notions:
\begin{lemma}\label{lemma:eqdgtensorfunctors}
Let $\T$ be a dg-category, two pseudo dg-tensor structures $\Gamma$, $\Lambda$ are tensor triangulated equivalent if there exists a pseudo dg-tensor functor $\Phi$ from $\Gamma$ to $\Lambda$ given by a dg-bimodule which is invertible under the tensor product of bimodules.
\begin{proof}
Suppose the bimodule $\Phi$ is invertible under the tensor product of dg-bimodules, so there exists $\Phi'$ such that $\Phi\otimes \Phi'$ is the identity $\T-Mod$ as a bimodule. \\ 
If $U$ and $U'$ are the units of $\Gamma$ and $\Lambda$ respectively, then we have equivalences $U\to \Phi\otimes\Phi' \otimes U$, and $\Gamma\otimes(\Phi\otimes\Phi')\otimes (\Phi\otimes \Phi')\to (\Phi\otimes\Phi')\otimes \Lambda$. Which are equivalent then to giving an equivalence $U\to U'$ and $\Gamma\to \Lambda$, and so for any $X,Y\in\T_{pe}$, $\Gamma\otimes X\otimes Y\simeq \Lambda\otimes X\otimes Y$, and similarly for the condition about the unit, and then $\hzero(\Gamma\otimes X\otimes Y) \cong \hzero(\Lambda\otimes X\otimes Y)$. \\
\end{proof}
\end{lemma}
As pointed out before, it is entirely possible that the structures differ at higher degrees and we only need the existence of pseudo dg-tensor functors $\Phi$ and $\Phi'$ such that $\hzero(\Phi\otimes \Phi')\simeq \hzero(\T-Mod)$. \\
Let us illustrate what we have so far with an example
\begin{exmp}\label{example:pseudodgtensorexceptional}
Let $X$ be a smooth projective variety and let us suppose that $\der$ has a full strong exceptional collection $\{E_{1},\dots,E_{m}\}$. In this case as we know, the object 
\[ E:=\bigoplus E_{i} \]
Is a compact generator and we have thus a homotopy equivalence of dg-categories 
\[ End(E)-Mod_{dg} \simeq D_{dg}^{b}(X) \] 
between the dg-category of dg-modules over $E$ and a dg-enhancement of $\der$.  \\
As the exceptional collection is strong there are no higher Ext groups and so this endomorphism algebra is supported in degree 0. Furthermore, we have the following description of this algebra
\[ \begin{pmatrix} M_{1,1} & M_{1,2} & \dots & M_{1,m} \\ 0 & M_{2,2} & \dots & M_{2,m} \\ \vdots & 0 & \ddots & \vdots \\ 0 & 0 & \dots & k \end{pmatrix}  \]
Where $M_{ii}=k$, and $M_{ij}$ is a right module over $M_{jl}$ for every $l$ and a left module over $M_{il}$ for any $l$. \\
As we can think of this dg-algebra concentrated in degree 0 as a regular k-algebra, we know that the category of dg-modules over it is simply the category of $E$-chain complexes. \\
The usual derived tensor product of $\der$ can be lifted to a 2-fold dg-bimodule over $E$ and it corresponds to the dg-bimodule given by 
\[ Hom(E\dtee E, E)\cong \bigoplus Hom(E_{i}\dtee E_{j}, E)\cong \bigoplus_{i,j,l} Hom(E_{i}\dtee E_{j}, E_{l}) \]
The best case scenario we can expect is for the variety to have a Picard group isomorphic to $\mathbb{Z}$ and the full strong exceptional collection to be composed of line bundles, in which case one might have a good chance of describing the 2-fold dg-bimodule $\Gamma$ corresponding to $\dtee$. \\

Let us consider an example of this situation and put $X=\mathbb{P}^{1}$, by using Beilinson's exceptional collection $\{\Ox_{X},\Ox_{X}(1)\}$. The 2-fold bimodule would then correspond to:
\[ \dots \to 0 \to \begin{pmatrix} 0 & 0 \\ 0 & 0 \\ 0 & 0 \\ k & 0 \end{pmatrix} \to \begin{pmatrix} k & k^{2} \\ 0 & k \\ 0 & k \\ 0 & 0 \end{pmatrix}\to 0\to \dots  \] 
Centered in degree 0.
Similarly we can calculate the unit object $U$ which  corresponds to
\[ \begin{pmatrix} k & k^2 \\ 0 & 0 \end{pmatrix}. \]
\end{exmp}
Following Hovey's paper we can deduce a few things about the classification of tensor triangulated categories on $\hzero(\ape)$ based on conditions about our dg-algebra $A$. \\
As an application of our encoding of tensor triangulated categories through perfect dg-tensor structures at a dg-enhancement we have the following result which is a derived version of a result of Hovey (\cite[Prop 4.1]{hovey2011additive} ).\\
\begin{defn}
Let $\Pp$ be a class of morphisms of chain complexes. We say that it is homotopically replete if for any $f\in \Pp$ such that there is a square
\[ \xymatrix{ X \ar[r]^{f} \ar[d] & Y \ar[d]  \\ X' \ar[r]^{g} & Y }  \]
Where the vertical morphisms induce homotopy equivalences, then the morphism $g$ is in $\Pp$
\end{defn}
\begin{prop}\label{prop:repleteclass}
Let $\T$ be a dg-category, and $\Gamma$ a perfect dg-pseudo tensor structure. Let $\Pp$ is a class of homotopically replete morphisms of chain complexes such that if $f\in \Pp$ then $X\otimes f$ and $f\otimes Y$, as morphisms of the underlying complexes, are in $\Pp$ for any pair of dg-modules $X,Y$. \\
If $f$ is a morphism of left $\T$-modules then $f\in \Pp$ if and only if the morphism $\Gamma\otimes f$ of the underlying complexes is in $\Pp$. Similarly $g$ a morphism of right $\T$-modules, then $g\in \Pp$ if and only if $g\otimes \Gamma$ is in $\Pp$.
\begin{proof}
We have that $g\simeq T-Mod\otimes g$ hence as the morphism $\Gamma\otimes U \otimes g \to \T-Mod\otimes g $ induces an isomorphism at the $\hzero$ -level, then $\Gamma\otimes U \otimes g$ is in $\Pp$. Similarly we have $\Gamma\otimes g \otimes U$. since $\Gamma\otimes g$ is in $\Pp$ by hypothesis, so is $g$. \\
The proof for the other structure on the left is similar. 
\end{proof}
\end{prop}
With this result in hand, it can be shown that
\begin{cor}\label{cor:homotopyreplete}
Let $\T$ be a dg-category and let $\Gamma$ be a perfect pseudo dg-tensor structure on $\T$. Then let $f$ be a morphism of $\T$-Modules and $g$ a morphism of $\T^{op}$-Modules. Then 
\begin{enumerate}
    \item $f\otimes \Gamma \simeq 0$ if and only if $f\simeq 0$. Similarly, $\Gamma\otimes g \simeq 0$ if and only if $g\simeq 0$
    \item $\hzero(f\otimes \Gamma)$ is an isomorphism if and only if $\hzero(f)$ is an isomorphism. Similarly for $g$.
    \item $\hzero(f\otimes \Gamma)$ is a surjection if and only if $\hzero(f)$ is a surjection.  
\end{enumerate}
\begin{proof}
Follows from morphisms $f\simeq 0$, homotopy equivalences and homotopy surjections forming a homotopically replete class. 
\end{proof}
\end{cor}
As a corollary to this then we have that 
\begin{cor}\label{cor:faithfulgamma}
Let $A$ be a dg-algebra and let $\Gamma$ be a perfect pseudo dg-tensor structure, then  $\Gamma$ is faithful as a $\hzero(A)$-module with either multiplication structure.
\begin{proof}
We let $a\in \hzero(A)$  a nonzero class. Then this induces a morphism $A\to A$ of dg-modules given by multiplication, which in turn induces a morphism $\Gamma\otimes \hzero(A)\to \Gamma\otimes \hzero(A)$, then the the fact that $a$ is not the zero class, implies the induced morphism is also not homotopically zero, thus $\Gamma$ is faithful as a $\hzero(A)$-module. 
\end{proof}
\end{cor}

%% file: DYCohomology.tex
\section{Davydov-Yetter cohomology}\label{3}
In \cite{davydov1997twisting, yetter1998braided} Davydov and Yetter introduced independently the concept of a deformation for tensor structures, one directly by deforming certain structural properties of the tensor category and the other by deforming monoidal functors and deforming the structural conditions it must admit. In this case the deformation of the identity functor seen as a monoidal functor takes the place of the deformation theory of the tensor structure directly. \\
Davydov-Yetter cohomology mainly parametrizes the associativity structural isomorphisms and provides an obstruction for these deformations. One important result that is obtained via this theory is Ocneanu's rigidity theorem which establishes that for fusion categories, which geometrically one can think of corresponding to a collection of points, admit no deformations of their tensor structures. \\
We quickly review the classical Davydov-Yetter cohomology construction for tensor categories and then we extend it to the case of categories of dg-modules over a dg-algebra $A$ equipped with a perfect pseudo dg-tensor structure. We will see that deformations of the associativity structural morphism can be described in an analogous manner to the abelian situation. \\
We finish the section with an investigation of the general deformation problem of tensor triangulated structures via our dg-enhancements and see how they relate to our extended Davydov-Yetter cohomology.  \\
For the remaining of this article and unless we say so, $k$ denotes an algebraically closed field of char 0.
\subsection{The abelian case}
Let us start by recalling a number of basic definition from the theory of tensor categories. We refer to the canonical reference \cite{etingof2016tensor}.
\begin{defn}
A tensor category is a k-linear rigid abelian monoidal category with biexact tensor product. \\
We say that it is finite if it is equivalent to a category of finite representations of a finite dimensional algebra.
\end{defn}
To a given finite tensor category we will assign a chain complex in the following way: \\
Through this section $\otimes$ denotes the monoidal product of a tensor category $\A$. \\
For any collection of objects $X_{1},\dots, X_{n}\in \A$ we denote by $\otimes^{n}(X_{1},\dots, X_{n})$ the full right parenthesization
\[ X_{1}\otimes (X_{2}\otimes ( \dots ( X_{n-1} \otimes X_{n})\dots ) \]
Similarly we denote by $^{n}\otimes$ the full left parenthesization
\[ ( \dots (X_{1}\otimes X_{2} )\dots )\otimes X_{n-1})\otimes X_{n} \]
For $n=1$ this assignment is simply the identity $\A\to \A$, and for $n=0$ then $\otimes^{0}$ is the constant functor determined by the monoidal unit in $\A$. \\
Given any parenthesization $X$ of a product of a collection of objects $X_{1},\dots, X_{n}$ there always exists a sequence of isomorphisms from $X$ to $\otimes^{n}(X_{1},\dots, X_{n})$ and to $^{n}\otimes(X_{1},\dots,X_{n})$. \\
Given any coherent endomorphism $f$ between any parenthesization of objects $X_{1},\dots,X_{n}$, we will write $\underline{f}$ to denote the morphism $f$ composed and precomposed by the coherent associative isomorphisms. \\
For example, if $f:X_{1}\otimes ((X_{2}\otimes X_{3})\otimes X_{4})\to (X_{1}\otimes X_{2})\otimes (X_{3}\otimes X_{4})$ then 
\begin{align*} \underline{f}:((X_{1}\otimes X_{2})\otimes X_{3})\otimes X_{4} \to X_{1}\otimes ((X_{2}\otimes X_{3})\otimes X_{4}) \\\to (X_{1}\otimes X_{2})\otimes (X_{3}\otimes X_{4}) \to  X_{1}\otimes (X_{2}\otimes (X_{3}\otimes X_{4})) \end{align*}\label{eq:diffdy}
This turns any such coherent morphism into a natural transformation in $Nat(^{n}\otimes,\otimes^{n})$. \\
The idea behind this operation is so that we can turn any morphism between different parenthesizations of a collection of objects $X_{1},\dots, X_{n}$ into the abelian group of morphisms $^{n}\otimes(X_{1},\dots, X_{n})\to \otimes^{n}(X_{1},\dots, X_{n})$. We however still need to keep track of the sign of this padding by associators, and so the morphism $\underline{f}$ comes with the signature $(-1)^{\mid sgn(f) \mid}$ where $sgn(f)$ denotes minimal number of associators $\alpha_{X,Y,Z}$ necessary to take a morphism $f$ into a morphism $f:^{n}\otimes(X_{1},\dots, X_{n})\to \otimes^{n}(X_{1},\dots, X_{n})$. For example the morphism \ref{eq:diffdy} has sign equal to 2 as we need to compose by the inverse of $Id_{X_{1}}\otimes \alpha_{X_{2},X_{3},X_{4}}$ and then by $\alpha_{X_{1}\otimes X_{2},X_{3},X_{4}}$. 
\begin{defn}\label{defn:davydovyettercomplex}
Let $\A$ be a finite tensor category, the Davydov-Yetter complex $DY^{\ast}(\A)$ is the chain complex defined in degree n by $Nat(^{n}\otimes,\otimes^{n})$. If $f\in DY^{n}$ a homogeneous element, the differential $d^{n}:DY^{n}\to DY^{n+1}$ is defined, on a set of objects $X_{1},\dots, X_{n+1}$ by 
\[ d^{n}(f):= Id_{X_{1}}\otimes f_{X_{1},\dots, X_{n}}+\Sigma_{i}(-1)^{i}f_{X{1},\dots,X_{i}\otimes X_{i+1},\dots, X_{n+1}} +(-1)^{n+1}f_{X_{1},\dots,X_{n}}\otimes Id_{X_{n+1}} \]
\end{defn}
It is a routine calculation to see that $d^{2}=0$ and so this forms a chain complex. \\
In degree 3 for example, we can calculate the  component of the Davydov-Yetter complex consists of natural transformations
\[ f_{X_{1},X_{2},X_{3}}:(X_{1}\otimes X_{2})\otimes X_{3} \to X_{1}\otimes (X_{2}\otimes X_{3}) \]
And has differential $d^{2}(f)_{X_{1},X_{2},X_{3},X_{4}}$  given by 
\[Id_{X_{1}}\otimes f_{X_{1},X_{2},X_{3}} - f_{X_{1}\otimes X_{2},X_{3},X_{4}} + f_{X_{1},X_{2}\otimes X_{3},X_{4}} - f_{X_{1},X_{2},X_{3}\otimes X_{4}} + f_{X_{1},X_{2},X_{3}}\otimes Id_{X_{4}} \]
With this definitions we can now define the Davydov-Yetter cohomology for finite tensor categories
\begin{defn}\label{defn:davydovyettercohomology}
Let $\A$ be a finite tensor category, the Davydov-Yetter cohomology $HDY^{\ast}(\A)$ is the cohomology of the the Davydov-Yetter complex $(DY^{\ast},d^{\ast})$
\end{defn}
Let us see by a hand calculation what the third cohomology group looks like
\begin{exmp}
The kernel of $d^{3}$ is composed of those natural transformations $f_{X_{1},X_{2},X_{3}}$ such that 
\begin{align*} Id_{X_{1}}\otimes f_{X_{1},X_{2},X_{3}} - f_{X_{1}\otimes X_{2},X_{3},X_{4}} + f_{X_{1},X_{2}\otimes X_{3},X_{4}} - f_{X_{1},X_{2},X_{3}\otimes X_{4}} \\ + f_{X_{1},X_{2},X_{3}}\otimes Id_{X_{4}}=0 \end{align*}
While the image of $d^{2}$ is
\[ Id_{X_{1}}\otimes f_{X_{1},X_{2}} - f_{X_{1}\otimes X_{2},X_{3}}+f_{X_{1},X_{2}\otimes X_{3}} - f_{X_{1},X_{2}}\otimes Id_{X_{3}} \]
\end{exmp}
As we said, we are interested in deformations of the associativity constraint of a finite tensor category $\A$, which is a coherent morphism $a_{X_{1},X_{2},X_{3}}:(X_{1}\otimes X_{2})\otimes X_{3} \to X_{1}\otimes (X_{2}\otimes X_{3})$. \\
We have then to say what we mean precisely by a deformation of this structure. In general the idea is that we should replace the coefficient ground field $k$ by an algebra that we interpret to represent a small neighborhood of it. Classically in the deformation theory of algebras we are interested in deforming by considering infinitesimal deformations of order n, so by replacing $k$ with the local ring $k[x]/x^{n+1}$. \\
We think of this ring as an augmented k-algebra, and so equipped with a morphism $k[x]/x^{n+1}\to k$ in the usual way. We always have a morphism in the other direction $i:k\to k[x]/x^{n+1}$.  \\
\begin{defn}\label{defn:deformationmonoidal}
An nth order deformation of a finite tensor category $\A$ over a field $k$ is a finite tensor category $\dfm{\A}$ over $k[x]/x^{n+1}$ such that there is a monoidal equivalence $\F$ to $\A$, after restriction of scalars $\_\otimes k$ of $Hom$ objects. \\
\end{defn}
\begin{defn}\label{defn:eqdeformations}
We say that two deformations are equivalent if there is a monoidal equivalence between them, with its underlying functor equal to the identity functor and such that the restriction of scalars of this monoidal equivalence is the identity after restriction of scalars. 
\end{defn}
In other words, as $Hom$ objects of the deformation $\dfm{\A}$ are modules over $k[x]/x^{n}$, the restriction of scalars functor given by tensoring by $k$ over $k[x]/x^{n}$ gives us an object of $k$-vector spaces. \\
To classify deformation classes of the associativity ( for example ) means that we are looking to classify associativity with coefficient in $k[x]/x^{n}$ which reduce to our original monoidal structure, and then classify the possible monoidal equivalences between them in the sense that we defined above. 
To say that we are deforming the associativity condition means that we are looking at deformations such that all the structural morphisms are given by extension of scalars except possibly for the associativity conditions. So for example deforming only the associativity condition means that the symmetry condition $X\otimes Y\to Y\otimes X$ in the deformation is exactly the extension by scalars given by $k[x]/x^n\to k$.
\begin{thm}\cite[Theorem 2.2]{crane1998deformations}\label{thm:davydovh3}
Let $\A$ be a finite tensor category, then there is a bijection between the 3rd Davydov-Yetter cohomology $HDY^{3}(\A)$ and equivalence classes of first order deformations of the associativity condition of $\A$.
\end{thm}
A fact is that the obstruction to these deformations to extend to infinitesimal deformations of second degree seem to be given by the 4th Davydov-Yetter cohomology group (cf. \cite{etingof2016tensor}). \\
When the category is semisimple we can say something about the deformation theory in this case:
\begin{thm}(Ocneanu's rigidity theorem)
Let $\A$ be a semisimple finite tensor category, then $HDY^{\ast}(\A)=0$.
\end{thm}
It is important to mention is that the hypothesis of the theorem cannot be omitted, as there exist explicit examples of non semisimple finite tensor categories with nontrivial Davydov-Yetter cohomology groups.\\
In \cite[Proposition 3.1]{davydov2018autoequivalences} it is shown that the collection of tensor autoequivalences of a finite tensor category is an affine algebraic group over $k$. \\
Concretely,
\begin{thm}\label{thm:algebraicgrouptensor}
Let $\A$ be a finite tensor category of $k$. Then $Aut(\A)$ has a natural structure of an affine algebraic group over $k$. 
\end{thm}
Our goal for the following is to define a Davydov-Yetter cohomology for tensor triangulated structures and study some of their properties in light of the work done in the abelian case. \\
\subsection{Deforming tensor triangulated structures}
In the previous section we saw how to construct a complex such that its cohomology behaves as the tangent space of a finite tensor category, in particular with respect to the associativity structural morphisms. Here we will construct a similar complex by using our perfect pseudo dg-tensor structures (see \ref{def:pseudodgtensor}). \\
From now on we will be interested in working with the category smooth proper derived noncommutative schemes (see \ref{defn:noncommuativescheme}) $NCSch^{pr}_{sm}$. As we have a homotopical equivalence between the derived categories of interest and a subcategory of such objects we exploit this to work our definitions at the dg-enhancement level. In this sense, whenever we refer to a dg-category $\T$ it is understood to be an object in $NSch^{pr}_{sm}$. We remark too that in general when working with dg-categories and dg-categories of dg-modules the change of basis operation $\_\dteee A$ for a dg-algebra $A$ does not in general preserve categories of perfect complexes, in the following however we try to keep notation as simple as possible and so we write $\otimes A$ for the operation $\widehat{\_\dteee A}_{pe}$ which takes a dg-category and sends it to the dg-category tensored by the dg-algebra $A$ (in the homotopy category of dg-categories  using the derived tensor product) after passing by its perfect closure. \\
Let us now commence to define our deformation complexes for our lifted tensor triangulated structures.
\begin{defn}
Let $\T$ be a dg-category and $\Gamma$ a perfect pseudo dg-tensor structure. For any $n\in \mathbb{N}$ we denote by $^{n}\Gamma$ the full left parenthesization of the 2-fold dg-bimodule $\Gamma$. This means we use the first multiplication to tensor $\Gamma$ with itself n-1 times by using this multiplication structure repeatedly. \\
This means that if $X_{1},\dots, X_{n}\in \T_{pe}$ then when taking tensor products we have
\[ ^{n}\Gamma:X_{1},\dots,X_{n}\mapsto \Gamma_{\Gamma,X_{1}}\otimes \Gamma_{\Gamma,X_{2}}\otimes \Gamma_{\Gamma,X_{3}}\otimes\dots \otimes \Gamma_{X_{n-1},X_{n}}\otimes X_{n}\otimes X_{n-1}  \]
Similarly we have $\Gamma^{n}$, which for a collection $X_{1},\dots, X_{n}\in \T_{pe}$
\[ ^{n}\Gamma:X_{1},\dots,X_{n} \mapsto \Gamma_{X_{1},\Gamma}\otimes X_{1}\otimes \Gamma_{X_{2},\Gamma}\otimes \dots \otimes \Gamma_{X_{n-1},X_{n}}\otimes X_{n-1}\otimes X_{n} \]
\end{defn}
\begin{defn}
Let $\Gamma$ be a perfect pseudo dg-tensor structure for the dg-algebra $A$ associated to a tensor triangulated structure on $\hzero(\ape)$. Let $DY_{dg}^{\ast,\ast}$ be the double chain complex defined by  
\[ DY_{dg}^{n,\ast}(A):= Bimod^{n}_{dg}(^{n-1}\Gamma,\;\Gamma^{n-1})^{\ast-1} \]
The complex of morphisms of dg-functors between the n-fold bimodule of totally left parenthesized product of $\Gamma$ with itself and the n-fold bimodule of totally right parenthesized product of $\Gamma$. \\
The vertical differential $d_{v}^{n,m}:DY_{dg}^{n,m}(A)\to DY_{dg}^{n,m+1}(A)$ is the differential given by the $Bimod(\;,\;)$ dg-functor. \\
The horizontal differential $d_{h}^{n,m}:DY_{dg}^{n,m}(A)\to DY_{dg}^{n+1,m}(A)$ is given, for $\eta^{\ast}\in DY_{dg}^{n,m}(A)$ and a collection $X_{1},\dots, X_{n}\in\ape$

\begin{align*}\label{eq:dy-diff} d_{h}^{n,m}(\eta)_{X_{1},\dots, X_{n}}= \\ Id_{X_{1}}\otimes \eta^{m}_{\Gamma \otimes_{\cc}X_{2}\otimes X_{3}, \dots, X_{n}} +\\  \Sigma_{i}(-1)^{i} \eta^{m}_{X_{1},\dots, \Gamma \otimes_{\cc}X_{i}\otimes X_{i+1},\dots, X_{n}}+(-1)^{n+1}\eta^{m}_{X_{1},\dots, X_{n-1}}\otimes Id_{X_{n}} \end{align*}


\begin{center}
\adjustbox{scale=0.7,center}{
\begin{tikzcd}
                 & \vdots                                                                                                       &  & \vdots                                                 &        \\
\cdots \arrow[r] & {DY_{dg}^{n,m+1}(A)=Bimod^{n}_{dg}(^{n-1}\Gamma,\;\Gamma^{n-1})^{m}} \arrow[u] \arrow[rr, "{d_{h}^{n,m+1}}"]                                                      &  & {DY_{dg}^{n+1,m+1}(A)} \arrow[u] \arrow[r]                & \cdots \\
\cdots \arrow[r] & {DY_{dg}^{n,m}(A)=Bimod^{n}_{dg}(^{n-1}\Gamma,\;\Gamma^{n-1})^{m-1}} \arrow[rr, "{d_{h}^{n,m}}"] \arrow[u, "{d_{v}^{n,m}}"'] &  & {DY_{dg}^{n+1,m}(A)} \arrow[u, "{d_{v}^{m+1,n}}"] \arrow[r] & \cdots \\
                 & \vdots \arrow[u]                                                                                             &  & \vdots \arrow[u]                                       &       
\end{tikzcd}
}
\end{center}
\end{defn}
A low degree example would be
\begin{center}
\begin{tikzcd}
                 & \vdots                                                                                                       &  & \vdots                                                 &        \\
\cdots \arrow[r] & {Bimod^{2}_{dg}(\Gamma,\Gamma)^{2}} \arrow[u] \arrow[rr, "{d_{h}^{2,3}}"]                                                      &  & {Bimod^{3}(^{2}\Gamma,\Gamma^{2})^{2}} \arrow[u] \arrow[r]                & \cdots \\
\cdots \arrow[r] & {Bimod^{2}(\Gamma,\Gamma)^{1}} \arrow[rr, "{d_{h}^{2,2}}"] \arrow[u, "{d_{v}^{2,1}}"'] &  & {Bimod^{3}(^{2}\Gamma,\Gamma^{2})^{1} } \arrow[u, "{d_{v}^{2,1}}"] \arrow[r] & \cdots \\
                 & \vdots \arrow[u]                                                                                             &  & \vdots \arrow[u]                                       &       
\end{tikzcd}
\end{center}

To make sense of the expression given by the differential, we proceed as in the abelian case, where the sacrifice that had to be made by working with non strict monoidal categories was that we had to introduce the padding construction to form additive groups of morphisms between different parenthesized products of a collection of objects. In this situation one could instead pass to a strict category using MacLane's coherence result and have that these objects are equal and there is no need for padding. \\
In our context however our only option is padding our morphisms, so let us explain how do we proceed: \\
Recall that the associativity coherence for a pseudo dg-tensor structure is given by a morphism  of dg-bimodules $\alpha_{X,Y,Z}:\Gamma_{\Gamma,Z}\otimes \Gamma_{X,Y}\to \Gamma_{X,\Gamma}\otimes \Gamma_{Y,Z}$ such that it becomes invertible when passing to $\hzero(\cc)$ for every triple of objects $X,Y,Z$. \\
Just as before whenever we have a morphism between two parenthesizations we would like to pad it so that it becomes a morphism in $Bimod^{\ast}(^{n}\Gamma,\Gamma^{n})(A)$. \\
Let us take then $n\in \mathbb{N}$ and a partition $n=n_{1}+\dots+n_{k}$. \\ 
A parenthesization of length $n\in \mathbb{N}$ of $\Gamma$ is a tensor product of factors of the form $^{n_{i}}\Gamma$ and $\Gamma^{n_{j}}$ along any of the two multiplications of $\Gamma$, for $0\leq i,j\leq k$. \\
Whenever we have a morphism between two such parenthesizations, \[f\in Bimod^{\ast}(^{n_{1}}\Gamma\otimes \dots \otimes \Gamma^{n_{k}}, \Gamma^{n'_{1}}\otimes\dots\otimes\Gamma^{n'_{l}})\], we would like to turn this into a morphism in $Bimod^{\ast}(^{n}\Gamma,\Gamma^{n})$ \\
The obstruction to do this as we did before is that we only know that the coherent morphism $\alpha$ is invertible when taking 0th chain cohomology $\hzero$. It is not enough to say that as this is invertible then in each step of the padding to consider a lift of the inverse $\alpha^{-1}$ as there are many and since we are only dealing with the truncation in low degree, the complexes in higher degree can change a lot from one another. So in reality what we must do is choose once and for all the inverse $\alpha^{-1}$ making $\alpha$ invertible. \\
Once the morphism is chosen we proceed as before, and then whenever we have for example a morphism
\[ f:\Gamma_{X_{1},\Gamma}\otimes\Gamma_{\Gamma,X_{4}}\otimes \Gamma_{X_{2},X_{3}}\to \Gamma_{\Gamma,\Gamma}\otimes\Gamma_{X_{1},X_{2}}\otimes\Gamma_{X_{3},X_{4}} \]
We can compose by products of the identity between dg-bimodules $\Gamma\to \Gamma$ and $\alpha$ and $\alpha^{-1}$. In this way we can obtain a morphism
\begin{align*} \underline{f} :\Gamma_{\Gamma, X_{4}}\otimes \Gamma_{\Gamma,X_{3}}\otimes\Gamma_{X_{1},X_{2}} \to \Gamma_{X_{1},\Gamma}\otimes\Gamma_{\Gamma,X_{4}}\otimes \Gamma_{X_{2},X_{3}} \\ \to \Gamma_{\Gamma,\Gamma}\otimes\Gamma{X_{1},X_{2}}\otimes\Gamma_{X_{3},X_{4}} \to \Gamma_{X_{1},\Gamma}\otimes\Gamma_{X_{2},\Gamma}\otimes  \Gamma_{X_{3},X_{4}}  \end{align*}
Now that we have defined our double complex, we can define
\begin{defn}\label{def:davydovyetterdg}
Let $\T$ be a dg-category and $\Gamma$ a perfect pseudo dg-tensor structure, the Davydov-Yetter complex is the total complex $Tot(DY_{dg}^{\ast,\ast})(\T)$ of the double complex $DY^{\ast,\ast}_{dg}(\T)$ with the usual differential $d_{tot}:=d_{v}+(-1)^{\mid v\mid \mid h\mid} d_{h}$. Where $\mid v \mid$ and $\mid h \mid $ denote the degree of the differentials $d_{h}$ and $d_{v}$.
\end{defn}
We put $^{0}\Gamma=\Gamma^{0}$ the bimodule $Id$ corresponding to the identity pseudo functor, and $^{-1}\Gamma=\Gamma^{-1}$ is the unit object $U$, and $^{n}\Gamma=\Gamma^{n}=0$ for any $n\leq -2$. \\
As before we calculate the cohomology of the dg-Davydov-Yetter complex
\begin{defn}\label{def:dgcohomologydavydovyetter}
Let $\T$ be a dg-category and $\Gamma$ a perfect pseudo dg-tensor structure and consider its Davydov-Yetter complex $Tot(DY_{dg}^{\ast,\ast})(\T)$. The total cohomology of this complex is the Davydov-Yetter cohomology of the perfect pseudo dg-tensor structure and we denote it by $HDY_{dg}^{\ast}(\T)$ or by $HDY_{dg}^{\ast}(\T,\Gamma)$ if there is ambiguity in which perfect pseudo dg-tensor structure is being considered.
\end{defn}
In lower degrees we have components of the total complex given by
\begin{align*} Tot^{3}(DY_{dg}^{\ast,\ast})=DY_{dg}^{3,0}(\T)\oplus DY_{dg}^{2,1}(\T)\oplus DY_{dg}^{1,2}(\T)\oplus DY_{dg}^{0,3}(\T) = \\
Bimod^{3}(^{2}\Gamma,\Gamma^{2})^{-1}\oplus Bimod^{2}(\Gamma,\Gamma)^{0}\oplus Bimod(\T_{pe},\T_{pe})^{1}\oplus \T_{pe}(U,U)^{2} \end{align*}
And 
\begin{align*} Tot^{4}(DY_{dg}^{\ast,\ast})= DY_{dg}^{4,0}(\T)\oplus DY_{dg}^{3,1}(\T)\oplus DY_{dg}^{1,3}(\T)\oplus DY_{dg}^{2,2}(\T)\\\oplus DY_{dg}^{3,1}(\T)\oplus DY_{dg}^{4,0} = \\
Bimod^{4}(^{3}\Gamma,\Gamma^{3})^{-1}\oplus Bimod^{3}(^{2}\Gamma,\Gamma^{2})^{0}\oplus Bimod^{2}(\Gamma, \Gamma)^{1} \\\oplus Bimod^{1}(\T_{pe},\T_{pe})^{2}\oplus \T_{pe}(U,U)^{3} \end{align*}
Let us calculate the action of the total differential in these degrees. If we let $\eta^{3}=(\eta^{3,1},\eta^{2,1},\eta{1,2},\eta^{0,3})$ then 
\begin{equation}\label{eq:diffd3}
\begin{split}
   d_{tot}^{3}(\eta^{3})=\\(d_{h}(\eta^{3,1}),d_{h}(\eta^{2,1})+d_{v}(\eta^{2,0}), d_{v}(\eta^{2,1}+d_{h}(\eta^{2,1}),\\d_{v}(\eta^{2,1})+d_{h}(\eta^{0,3}), d_{v}(\eta^{0,3}))= \\
   (\eta^{3,0}_{X,Y,Z}\otimes Id_{W} - \eta^{3,0}_{\Gamma\otimes X\otimes Y, Z,W}+\eta^{3,0}_{X,\Gamma\otimes y\otimes Z, W}-\eta^{3,0}_{X,Y,\Gamma \otimes Z\otimes W}+Id_{X}\otimes \eta^{3,0}_{Y,Z,W}, \\
   \eta^{2,1}_{X,Y}\otimes Id_{Z}-\eta^{2,1}_{\Gamma\otimes X\otimes Y, Z}+ \eta^{2,1}_{X,\Gamma\otimes Y\otimes Z}+ Id_{X}\otimes\eta^{2,1}_{Y,Z}+d_{v}(\eta^{3,0}),\\
   \eta^{1,2}_{X}\otimes Id_{Y}-\eta^{1,2}_{\Gamma\otimes X\otimes Y}+ Id_{X}\otimes \eta^{2,1}_{Y} + d_{v}(\eta^{2,1}), \\
   d_{v}(\eta^{0,3}))
   \end{split}
   \end{equation}

And to calculate the kernel of $d_{tot}^{4}$, we see that it corresponds to those $\eta^{4}=(\eta^{4,0},\eta^{3,1},\eta^{2,2},\eta^{1,3},\eta^{0,4})$ such that
\begin{equation}\label{eq:diffd4}
\begin{split}
    d_{tot}^{4}(\eta^{4})= \\(d_{h}(\eta^{4,0}),-d_{h}(\eta^{3,1})+d_{v}(\eta^{4,0}), d_{h}(\eta^{2,2})+d_{v}(\eta^{3,1}), d_{h}(\eta^{1,3})+d_{v}(\eta^{2,2}),\\d_{h}(\eta^{0,4})+d_{v}(\eta^{1,3},d_{v}(\eta^{0,4})) = \\
    (Id_{X}\otimes \eta^{4,0}_{Y,Z,W,R} - \eta^{4,0}_{\Gamma\otimes X\otimes Y,Z,W,R}+\eta^{4,0}_{X,\Gamma\otimes Y\otimes Z,W,R}-\eta^{4,0}_{X,Y,\Gamma\otimes Z\otimes W, R}\\+ \eta^{4,0}_{X,Y,Z,\Gamma\otimes W\otimes R} - \eta^{4,0}_{X,Y,Z,W}\otimes Id_{R}, \\ -Id_{X}\otimes \eta^{3,1}_{Y,Z,W}+\eta^{3,1}_{\Gamma\otimes X\otimes Y, Z,W} - \eta^{3,1}_{X,\Gamma\otimes Y\otimes Z, W} +\eta^{3,1}_{X,Y,\Gamma\otimes Z\otimes W}\\ -\eta^{3,1}_{X,Y,Z}\otimes Id_{W} + d_{v}(\eta^{4,0}), \\ Id_{X}\otimes\eta^{2,2}_{Y,Z} - \eta^{2,2}_{\Gamma\otimes X\otimes Y, Z} + \eta^{2,2}_{X,\Gamma\otimes Y\otimes Z } - \eta^{2,2}_{X,Y}\otimes Id_{Z} + d_{v}(\eta^{3,1}), \\
    Id_{X}\otimes \eta^{1,3}_{Y} - \eta^{1,3}_{\Gamma\otimes X\otimes Y}+ \eta^{1,3}_{X}\otimes Id_{Y} +d_{v}(\eta^{2,2}), \\
    Id_{X}\otimes U + d_{v}(\eta^{1,3}), d_{v}(\eta^{0,4}))=0
\end{split}
\end{equation}
We should describe what we mean by a deformation of the structure,  $\Gamma$. For this we need 
\begin{defn}\label{def:deformationdg}
Let $\T$ be a dg-category, $\Gamma$ a perfect pseudo dg-tensor structure on $\T$. An n-th order deformation of $\Gamma$ consists of a perfect pseudo dg-tensor structure $\dfm{\Gamma}$ on $\T\otimes_{k}k[x]/x^{n+1}$ such that \[i^{\ast}\dfm{{\Gamma}}:=\dfm{\Gamma}\otimes_{k[x]/x^{n+1}}k\] is a perfect pseudo dg-tensor structure equivalent to $\Gamma$. 
\end{defn}
\begin{defn}\label{def:eqdgdeformations}
We say two n-th order infinitesimal deformations of a perfect pseudo dg-tensor structure $\Gamma$ are equivalent if there is a pseudo dg-tensor functor $\Phi$ in $\T\otimes_{k} k[x]/x^{n}$ such that its restriction $\Phi\otimes_{k[x]/x^{n+1}}k$ is equivalent to the pseudo dg-tensor functor given by the identity dg-bimodule $\T_{pe}$.
\end{defn}
Just as before, we will say that a deformation of a perfect pseudo dg-structure $\Gamma$ is a deformation of the associativity condition if the deformation $\dfm{\Gamma}$ has as structure morphisms for the unit and symmetry conditions equivalent equivalent to $\ell\otimes k[x]/x^{n+1}$,$r\otimes k[x]/x^{n+1}$, and $c\otimes k[x]/x^{n+1}$, while we allow for the associativity coherence condition to possibly be different.  \\
The following is an analogue result to Theorem \ref{thm:davydovh3}
\begin{thm}\label{thm:deformationdg}
Let $\T$ be a dg-category and let $\Gamma$ be a pseudo dg-tensor structure on $\T$. Then to any element of  $HDY^{4}_{dg}(\T)$ we can associate an  equivalence class of infinitesimal deformations of order 1 of the associativity condition of $\Gamma$.
\begin{proof}
Let us recall that 
\begin{align*} DY^{4}_{dg}(\T) = Bimod^{4}(^{3}\Gamma,\Gamma^{3})^{-1}\oplus Bimod^{3}(^{2}\Gamma,\Gamma^{2})^{0} \\ \oplus Bimod^{2}(\Gamma,\Gamma)^{1}\oplus Bimod^{1}(\T_{pe},\T_{pe})^{2}\oplus \T_{pe}(U,U)^{3} \end{align*}
As we calculated, the kernel of $d^{4}_{dg}$ consists of those $\eta\in DY^{4}_{dg}(\T)$ such that the equation \ref{eq:diffd4} is equal to zero. \\
We see that the first component
\[ Id_{X}\otimes \eta^{4,0}_{Y,Z,W,R} - \eta^{4,0}_{\Gamma\otimes X\otimes Y,Z,W,R}+\eta^{4,0}_{X,\Gamma\otimes Y\otimes Z,W,R}-\eta^{4,0}_{X,Y,\Gamma\otimes Z\otimes W, R}\\+ \eta^{4,0}_{X,Y,Z,\Gamma\otimes W\otimes R} - \eta^{4,0}_{X,Y,Z,W}\otimes Id_{R} \]
is a morphism of 4-fold bimodules, $\eta^{4,0}$ satisfying an hexagon condition. \\
The second component
\[ -Id_{X}\otimes \eta^{3,1}_{Y,Z,W}+\eta^{3,1}_{\Gamma\otimes X\otimes Y, Z,W} - \eta^{3,1}_{X,\Gamma\otimes Y\otimes Z, W} +\eta^{3,1}_{X,Y,\Gamma\otimes Z\otimes W}\\ -\eta^{3,1}_{X,Y,Z}\otimes Id_{W} + d_{v}(\eta^{4,0})
\] 
on the other hand, is a 3-fold bimodule morphism, $\eta^{3,1}\in Bimod^{3}(^{2}\Gamma,\Gamma^{2})^{0}$ satisfying the pentagon diagram condition up to the homotopy $d_{v}(\eta^{4,0})$. This is precisely the condition we require as a coherence condition for the associator morphism of a pseudo dg-tensor structure. \\
To be more precise, we will consider a perfect pseudo dg-tensor structure on $\T\otimes k[x]/x^{2}$ given by the 2-fold dg-bimodule $\Gamma^{k[e]}:= \Gamma\otimes_{k}k[x]/x^{2}$ where every one of the structural morphisms of Definition \ref{def:pseudodgtensor} are given by the extension of scalars $\otimes_{k}k[x]/x^{2}$ except the associator. \\
Our goal is to define a new associator $\dfm{\alpha}$ which will restrict back to the associator $\alpha$. \\
Let us write then $\dfm{\alpha}:= \alpha+ \eta^{3,1}x\in Bimod(^{2}\Gamma^{k[e]},\Gamma^{k[3]}\,^{(2)})^{0}$. And so, obtain a natural morphism
\[ \dfm{\alpha}_{X,Y,Z}:\Gamma^{k[e]}_{\Gamma^{k[e]},Z}\otimes \Gamma^{k[e]}_{X,Y}\to \Gamma^{k[e]}_{X,\Gamma^{k[e]}}\otimes \Gamma^{k[e]}_{Y,Z}. \]
We need to check that this morphism satisfies the pentagon identity.
\[ -Id_{X}\otimes \dfm{\alpha}_{Y,Z,W}\circ \dfm{\alpha}_{X,\Gamma^{k[e]}_{Y,Z},W}\circ \dfm{\alpha}_{X,Y,Z}\otimes Id_{W} + \dfm{\alpha}_{X,Y,\Gamma^{k[e]}_{Z,W}}\circ\dfm{\alpha}_{\Gamma^{k[e]}_{X,Y},Z,W} = d_{v}(\eta^{4,0}) \]
Recall to add morphisms between different parenthesizations we need first to pad the morphism in the sense that we need to compose and pre-compose by the associativity morphism $\alpha\otimes k[x]/x^{2}$ and a fixed choice of an inverse $\alpha^{-1}\otimes k[x]/x^{2}$. We do this in such a way that addition of morphisms is always between the leftmost parenthesization and the rightmost one. \\
We proceed as in the abelian case and see that since the associator $\alpha$ of the pseudo dg-tensor structure $\Gamma$ already satisfies the pentagon up to homotopy and the relation $x^{2}$, our associator $\dfm{\alpha}$ satisfies the pentagon diagram up to homotopy because $\alpha$ already satisfies the condition, and as we are working with coefficients in $k[x]/x^{2}$ then all that is left is the expression involving $\alpha$ and the morphisms $\eta^{3,1}$ which corresponds to our padding of morphisms and then we obtain the condition $d_{h}^{3}(\eta^{3,1})=d_{v}(\eta^{4,0})$. \\
Now we would like to see that two such deformations $\dfm{\alpha}, \dfm{\alpha}'$ given as above are equivalent then they come from a pseudo dg-tensor functor induced equivalence of perfect pseudo dg-tensor structures. 
As in the abelian case the underlying functor we are looking for is the identity functor and so our dg-bimodule $\Phi$ is nothing but the bimodule $\T_{pe}$. \\
This means that we are looking for a morphism 
\[\Gamma\otimes k[x]/x^{2}\simeq \Gamma\otimes\T_{pe}\otimes\T_{pe}\otimes k[x]/x^{2}\to \T_{pe}\otimes \Gamma\otimes k[x]/x^{2}\simeq \Gamma\otimes k[x]/x^{2}\]
of the form 
\[ Id+\beta x\]
where 
\[\beta: \Gamma\to \Gamma\]
is a morphism of dg-bimodules. \\
We write $\beta$ then as any dg-bimodule morphism such that $\beta^{0}=\eta^{2,1}$, and writing the associativity condition for the pseudo dg-tensor we see that $\dfm{f}:=Id+\beta x$ satisfies this condition if $\eta^{2,1}$ is in the image of $d_{tot}^{2}$, as the second component of $d_{tot}^{2}$ gives us precisely the coherence diagram up to homotopy that $\dfm{f}$ has to satisfy. \\
Indeed as the identity dg-bimodule $\T_{pe}$ is a pseudo dg-tensor functor the identity $\T_{pe}\to \T_{pe}$ satisfies the associativity condition and as in the abelian case we see that the only remaining morphisms are those composed with the associativity condition of the pseudo dg-tensor structure, this in turn corresponds to a padding operation and thus we obtain precisely the expression in the second component of the image of $d_{tot}^{2}$.

\end{proof}
\end{thm}
 
The converse of this theorem don't seem to hold in general and likely require either a more general and coherent setting in which the deformation and space of tensor structures takes place, or under stricter conditions for the pseudo dg-tensor structure itself. At this point we ignore what the higher coherence conditions appearing both in the kernel and the image of the differentials $d_{tot}$ represent in the context of deformations of the lifts of tensor triangulated structures. In all likeness an approach where the lift is meant to produce a tensor structure itself in the dg-enhancement is the correct setting in which one ought to take these deformations. Having said that as our motivation was kept in line with the tensor structure at the triangulated category level and this is the reason for the brute truncation of these tensor structures. \\
In the abelian case we saw a few more things, namely we had Ocneanu's rigidity theorem, or \cite[Proposition 3.21]{batanin2020cosimplicial} exhibiting Davydov-Yetter cohomology as a tangent space of the moduli functor of tensor structures on a given finite tensor category $\A$. Additionally in the same paper Batanin and Davydov investigate a good deal of the deformation theory of tensor structures, like a Lie algebra structure with many of the usual expected properties. We refer to their work for further details. \\
In our case we expect that a dg version of Ocneanu's rigidity theorem for our perfect pseudo dg-tensor category structures should hold in the same spirit, but a concrete formulation of such a theorem must then deal again with the precise nature of the truncation we are performing in our structure. \\
Let us mention the work of Panero-Shoikhet in \cite{panero2022category} in which they do a thorough analysis of the deformation theory of monoidal dg-categories through a Davydov-Yetter cohomology construction. In contrast with our approach, theirs deals with entirely with the strict theory of dg-categories, that is, without considering homotopy equivalences. \\
The concept of monoidal dg-category seems to be relatively elusive in the literature, with different authors meaning different things, often but not always with an implicit understanding that these definitions should be equivalent in some way. \\
Our pretriangulated dg-categories equipped with a given pseudo dg-tensor category structure can be, as hinted already a number of times through this work, as a truncated version of a true homotopy monoidal dg-category in the Morita model category structure. The following formal statement reflecting this was suggested to the author by Bertrand Toën. 
\begin{thm}\label{thm:truncationmonoidal}
Let $Ho_{2}(dg-cat_{k})$ be the 2-category given by the homotopy category of the 2-truncation $\tau_{\leq 2}N(dg-cat_{k})[W^{-1}]$ of the $(\infty,1)$-category $N(dg-cat_{k})[W^{-1}]$ of dg-categories. Let $\M$ be a dg-category equivalent to a derived noncommutative scheme then any perfect pseudo dg-tensor structure $\Gamma$ on $\M$ induces an associative monoid in $Ho_{2}(dg-cat_{k})$ and to any associative monoid structure on $\M$ in $Ho_{2}(dg-cat_{k})$ induces a possibly non unique perfect pseudo dg-tensor structure on $\M$.
\end{thm}
Where $W$ is the class of Morita equivalences (See \cite{tabuada2005invariants} for the definition of the Morita model category structure). \\
\\
A moduli space interpretation of this deformation theory would be of most interest to us and we hope to come back to the question in future work. On the other hand, we can give a dg version of Theorem \ref{thm:algebraicgrouptensor} which can be thought as a local moduli space result. \\
\begin{thm}\label{thm:representationttsgamma}
Let $\T$ be a dg-category, $A$ a k-algebra, a perfect module $U$ and $\Gamma$ a 2-fold dg-bimodule over $\T\otimes^{\mathbb{L}} A$ with finite global dimension. Then the set $TTS_{A}(\Gamma)$ of perfect pseudo dg-tensor structures over $\T\otimes^{\mathbb{L}} A$ which has $\Gamma$ as a 2-fold dg-bimodule and $U$ as a unit, has a structure of a quotient of an affine scheme by an algebraic affine scheme.
\begin{proof}
As both $\Gamma$ and $U$ are fixed, what we are looking for is simply families of morphisms $\alpha:\Gamma\otimes\Gamma\to \Gamma\otimes\Gamma$, $u:\Gamma\otimes U\to \T_{pe}$, and $c:\Gamma\to \Gamma$ we use that we are under the assumption that our dg-category $\T$ is proper and so locally perfect, which means that every complex of morphisms is bounded and of finite dimension as $k$-vector spaces. \\
The polynomial equations that determine morphisms $\alpha$, $c$ and $u$ together with the differentials in every chain complex determine then an affine scheme and so do they when we restrict them by imposing the coherence conditions. \\
For two such structures, $\Gamma_{1}$, $\Gamma_{2}$ to be equivalent, we need an invertible pseudo dg-tensor functor $\Phi$ from, say, $\boxtimes_{1}$ to $\boxtimes_{2}$ ( using Lemma \ref{lemma:eqdgtensorfunctors} ). Meaning that we have $\Phi'$ from $\boxtimes_{2}$ to $\boxtimes_{1}$ such that $\Phi\otimes \Phi'$ is isomorphic to the dg-bimodule $\T_{pe}$. \\
As the dg-bimodule $\Gamma$ is the same in both structures $\Gamma_{1}$ and $\Gamma_{2}$, any such pseudo dg-tensor functor is determined by its structural morphisms
\[ \Gamma\to \Gamma \]
and
\[ \Gamma\otimes U\to \Gamma \]
Just as before the complexes are bounded and finite dimensional $k$-vector spaces and then they determine an affine scheme. Using composition as the group operation, taking the identity pseudo dg-tensor functor $T_{pe}$ as the identity element and noting that every pseudo dg-tensor functor has an inverse, we get that the affine scheme of these equivalences forms an algebraic group. \\
We have then a quotient of the space of pseudo dg-tensor structures with 2-fold dg-bimodule $\Gamma$ and unit $U$ by the affine group scheme of equivalences between them. 
\end{proof}
\end{thm}